\documentclass[fleqn,12pt]{article}
\usepackage{makeidx}
\usepackage{graphicx}
\usepackage{multicol}
\usepackage[bottom]{footmisc}
\usepackage{amssymb}
\usepackage{graphicx}
\usepackage{amsmath}
\usepackage{url}
\usepackage{algorithmicx}
\usepackage{times}
\begin{document}

\title{Recursive computation of spherical harmonic rotation coefficients of
large degree}
\author{Nail A. Gumerov$^{1}$ and Ramani Duraiswami$^{1, 2}$ \\
$^1$ Institute for Advanced Computer Studies\\$^2$ Department of Computer Science, \\ University of Maryland,
College Park. \\ \url{gumerov@umiacs.umd.edu},  \url{http://www.umiacs.umd.edu/~gumerov} \\ \url{ramani@umiacs.umd.edu},  \url{http://www.umiacs.umd.edu/~ramani}
}
\date{29 March 2014} 
\maketitle
\begin{abstract}
Computation of the spherical harmonic rotation coefficients or elements of
Wigner's d-matrix is important in a number of quantum mechanics and
mathematical physics applications. Particularly, this is important for the
Fast Multipole Methods in three dimensions for the Helmholtz, Laplace  and related equations, if rotation-based decomposition of translation operators are used. In these and related problems related to representation of functions on a sphere via spherical harmonic expansions, computation of the rotation
coefficients of large degree $n$ (of the order of thousands and more) may be necessary. Existing algorithms for their computation, based on recursions,  are usually unstable, and do not extend to $n$. We develop a new recursion and study its behavior for large degrees, via computational and asymptotic analyses. Stability of this recursion was studied based on a novel application of the Courant-Friedrichs-Lewy condition and the von Neumann method for stability of finite-difference schemes for solution of PDEs. A recursive algorithm of minimal complexity $O\left( n^{2}\right)$ for degree $n$ and FFT-based algorithms of complexity $O\left( n^{2}\log n\right) $ suitable for computation of rotation coefficients of large degrees are proposed, studied numerically, and cross-validated. It is shown that the latter algorithm can be used for $n\lesssim 10^{3}$ in  double precision, while the former algorithm was tested for large $n$ (up to $10^{4}$ in our experiments) and demonstrated better performance and accuracy compared to the FFT-based algorithm. 


\vspace*{0.5in}

{University of Maryland Institute for Advanced Compute Studies Technical Report, UMIACS-TR-2014-04 \\ \\ University of Maryland Department of Computer Science Technical Report,  CS-TR-5037}
\end{abstract}
\newpage
\tableofcontents
\newpage
\section{Introduction}
Spherical harmonics form an orthogonal basis for the space of square integrable functions defined over the unit
sphere, $S_{u}$, and have important application for a number of problems of
mathematical physics, interpolation, approximation, and Fourier analysis on
the sphere. Particularly, they are eigenfunctions of the Beltrami operator
on the sphere, and play a key role in the solution of Laplace, Helmholtz, and
related equations (polyharmonic, Stokes, Maxwell, Schroedinger, etc.) in
spherical coordinates. Expansions of solutions of these equations via
spherical basis functions, whose angular part are the spherical
harmonics (multipole and local expansions), are important in the fast multipole
methods (FMM) \cite{Greengard1987:JCP, Rokhlin93:ACHA, Epton95:SISC, Chew2001, Gumerov05:Book}.

The FMM for the Helmholtz equation as well as other applications in geostatistics require
operations with expansions which involve large numerical values of the
maximum degree of the expansion, $p$. This may reach several thousands, and
the expansions, which have $p^{2}$ terms, will have  millions of coefficients. For the FMM for the Helmholtz equation, such  expansions arise when the  domain has size of the order of $M\sim 100$
wavelengths, and for convergence we need $p \sim O(M)$. These large expansions need to be  translated (change of the origin of
the reference frame). Translation operators for truncated expansions of
degree $p$ ($p^{2}$ terms) can be represented by dense matrices of size $%
\left( p^{2}\right)^{2}=p^{4}$ and, the translation can be performed via matrix vector product with cost $O\left( p^{4}\right)$.
Decomposition of the translation operators into rotation and coaxial
translation parts (the RCR-decomposition: rotation-coaxial
translation-back rotation) \cite{White96:JChemPhys, Gumerov05:Book} reduces
this cost to $O\left( p^{3}\right)$.
While translation of expansions can be done with asymptotic complexity $%
O\left( p^{2}\right) $ using diagonal forms of the translation operators 
\cite{Rokhlin93:ACHA}, use of such forms  in the
multilevel FMM requires additional operations, namely, interpolation and
anterpolation, or filtering of spherical harmonic expansions, which can be
performed for $O\left( p^{2}\log p\right) $ operations. The practical
complexity of such filtering has large asymptotic constants, so that $%
\left( p^{3}\right)$ methods are competitive with asymptotically faster
methods for $p$ up to several hundreds \cite{Suda02:MathComp}. So, wideband
FMM for the Helmholtz equation for such $p$ can be realized in different
ways, including \cite{Cheng06:JCP} formally scaled as $O\left( p^{2}\log
p\right) $ and \cite{Gumerov09:JASA}, formally scaled as\ $O\left(
p^{3}\right)$, but with comparable or better performance for the asymptotically slower method.

Rotation of spherical harmonic expansions is needed in several
other applications (e.g., \cite{park09}, and is interesting from a mathematical point of view,
and has deep links with group theory \cite{Vilenkin68:Book}. Formally, expansion of degree $p$ can
be rotated for the expense of $p^{3}$ operations, and, in fact, there is a
constructive proof that this can be done for the expense of $O\left(
p^{2}\log p\right) $ operations \cite{Gumerov05:Book}. The latter is related
to the fact that the rotation operator (matrix) can be decomposed into the
product of diagonal and Toeplitz/Hankel matrices, where the matrix-vector
multiplication involving the Toeplitz/Hankel matrices can be performed for $%
O\left( p^{2}\log p\right) $ operations using the FFT. There are two issues
which cause difficulty with the practical realization of such an algorithm.
First, the matrix-vector products should be done for $O(p)$ matrices of
sizes $O\left( p\times p\right) $ each, so for $p\sim 10^{2}-10^{3}$ the
efficiency of the Toeplitz matrix-vector multiplication of formal complexity 
$O(p\log p)$ per matrix involving two FFTs is not so great compared to a
direct matrix-vector product, and so the practical complexity is comparable
with $O\left( p^{2}\right) $ brute-force multiplication due to large enough
asymptotic constant for the FFT. Second, the decomposition shows poor
scaling of the Toeplitz matrix (similar to Pascal matrices), for which
renormalization can be done for some range of $p$, but is also algorithmically
costly \cite{Tang04:TechRep}.

Hence, from a practical point of view $O\left( p^{3}\right) $ methods of
rotation of expansions are of interest. Efficient $O\left( p^{3}\right) $
methods are usually based on direct application of the rotation matrix to
each rotationally invariant subspace, where the the rotation coefficients
are computed via recurrence relations. There are numerous recursions, which
can be used for computation of the rotation coefficients (e.g. \cite{Biedenharn81:Book, Ivanic96:JPhysChem, Gumerov03:SISC, Dachsel06:JChemPhys}, see also the review in \cite{Gimbutas09:JCP}), and
some of them were successfully applied for solution of problems with
relatively small $p$ ($p\lesssim 100$). However, attempts to compute
rotation coefficients for large $p$ using these recursions face numerical instabilities. An $O\left( p^{3}\right) $ method
for rotation of spherical expansion, based on pseudospectral projection,
which does not involve explicit computation of the rotation coefficients was
proposed and tested in \cite{Gimbutas09:JCP}. The rotation coefficient values are however needed in some applications. For example for the finite set of fixed angle rotations encountered in the FMM, the rotation
coefficients can be precomputed and stored. In this case the algorithm
which simply uses the precomputed rotation
coefficients is faster than the  method proposed in \cite%
{Gimbutas09:JCP}, since brute force matrix-vector multiplications do not
require additional overheads related to spherical harmonic evaluations and
Fourier transforms \cite{Gimbutas09:JCP}, and are well optimized on hardware.

We note that almost all studies related to computation of
the rotation coefficients that advertise themselves as ``fast and stable'', in
fact, do not provide an actual stability analysis. ``Stability'' then is 
rather a reflection of the results of numerical experiments conducted for
some limited range of degrees $n$. Strictly speaking, all  algorithms
that we are aware of for this problem are not proven stable in the strict sense -- that the error in computations is not
increasing with increasing $n$. While there are certainly unstable schemes,
which ``blow up'' due to exponential error growth, there are some unstable
schemes with slow error growth rate $O\left( n^{\alpha }\right) $ at
large $n$. 

In the present study we investigate the behavior of the rotation
coefficients of large degree and propose a ``fast and stable'' $O\left(
p^{3}\right) $ recursive method for their computation, which numerically is
much more stable than other algorithms based on recursions used in the
previous studies. We found regions where the recursive processes used
in the present scheme are unstable as they violate a
Courant-Friedrichs-Lewy (CFL) stability condition \cite{Courant28:MathAn}.
The proposed algorithm manages this. We also show that in the
regions which satisfy the CFL condition the recursive computations despite
being formally unstable have a slow error growth rate. Such conclusion comes
partially from the well-known von Neumann stability analysis \cite%
{Charney50:QJG} combined with the analysis of linear one-dimensional
recursions and partially from the numerical experiments on noise
amplification when using the recursive algorithm. The proposed algorithm was
tested for computation of rotation coefficients of degrees up to $%
n=10^{4}$, without substantial constraints preventing their use for larger $n$. We also proposed
and tested a non-recursive FFT-based algorithm of complexity $O\left(
p^{3}\log p\right) $, which despite larger complexity and higher errors
than the recursive algorithm, shows good results for $n\lesssim 10^{3}$ and can be used for validation (as we did) and other purposes.
\section{Preliminaries}
\subsection{Spherical harmonic expansion}
Cartesian coordinates of points on the unit sphere are related to the angles
of spherical coordinates as 
\begin{equation}
\mathbf{s}=\left( x,y,z\right) =\left( \sin \theta \cos \varphi ,\sin \theta
\sin \varphi ,\cos \theta \right) ,  \label{1}
\end{equation}
We consider functions $f\in L_{2}\left( S_{u}\right) $ of bandwidth $p$,
means expansion of $f\left( \theta ,\varphi \right) $ over spherical
harmonic basis can be written in the form 
\begin{equation}
f\left( \theta ,\varphi \right)
=\sum_{n=0}^{p-1}\sum_{m=-n}^{n}C_{n}^{m}Y_{n}^{m}\left( \theta ,\varphi
\right) ,  \label{2}
\end{equation}
where orthonormal spherical harmonics of degree $n$ and order $m$ are
defined as 
\begin{align}
Y_{n}^{m}\left( \theta ,\varphi \right) & =(-1)^{m}\sqrt{\frac{2n+1}{4\pi }%
\frac{(n-\left| m\right| )!}{(n+\left| m\right| )!}}P_{n}^{\left| m\right|
}(\cos \theta )e^{im\varphi },  \label{3} \\
n& =0,1,2,...;\qquad m=-n,...,n.  \notag
\end{align}
Here $P_{n}^{m}\left( \mu \right) $ are the associated Legendre functions,
which are related to the Legendre polynomials, and can be defined by the
Rodrigues' formula
\begin{eqnarray}
P_{n}^{m}\left( \mu \right) &=&\left( -1\right) ^{m}\left( 1-\mu ^{2}\right)
^{m/2}\frac{d^{m}}{d\mu ^{m}}P_{n}\left( \mu \right) ,\quad
n=0,1,2,...,\quad m=0,1,2,...  \label{4} \\
P_{n}\left( \mu \right) &=&\frac{1}{2^{n}n!}\frac{d^{n}}{d\mu ^{n}}\left(
\mu ^{2}-1\right) ^{n},\quad n=0,1,2,...  \notag
\end{eqnarray}
The banwidth $p$ can be arbitrary, and the fact
that we consider finite number of harmonics relates only to computations. We note also that different
authors use slightly different definitions and normalizations for the spherical harmonics. Our definition is consistent with that of \cite{Epton95:SISC}. Some discussion on
definitions of spherical harmonics and their impact on translation relations
can be also found there. Particularly, for spherical
harmonics defined as 
\begin{align}
\widetilde{Y}_{n}^{m}\left( \theta ,\varphi \right) & =\sqrt{\frac{2n+1}{%
4\pi }\frac{(n-m)!}{(n+m)!}}P_{n}^{m}(\cos \theta )e^{im\varphi },
\label{4.1} \\
n& =0,1,2,...;\qquad m=-n,...,n,  \notag
\end{align}
\begin{equation}
Y_{n}^{m}\left( \theta ,\varphi \right) =\epsilon _{m}\widetilde{Y}%
_{n}^{m}\left( \theta ,\varphi \right) ,  \label{4.2}
\end{equation}
where 
\begin{equation}
\epsilon _{m}=\left\{ 
\begin{array}{c}
(-1)^{m},\quad m\geqslant 0, \\ 
1,\quad m<0.%
\end{array}
\right.  \label{r6}
\end{equation}
Hence one can expect appearance of factors $\epsilon _{m}$ in relations used
by different authors.
\subsection{Rotations}
There are two points of view on rotations, active (alibi), where vectors are
rotated in a fixed reference frame, and passive (alias), where vectors are
invariant objects, but the reference frame rotates and so the coordinates of
vectors change. In the present paper we use the latter point of view,
while it is not difficult to map the relations to the active view
by replacing rotation matrices  by their transposes
(or inverses). 

\begin{figure}[htb]
\center  \includegraphics[width=0.9\textwidth, trim=0 1.4in 0 0]{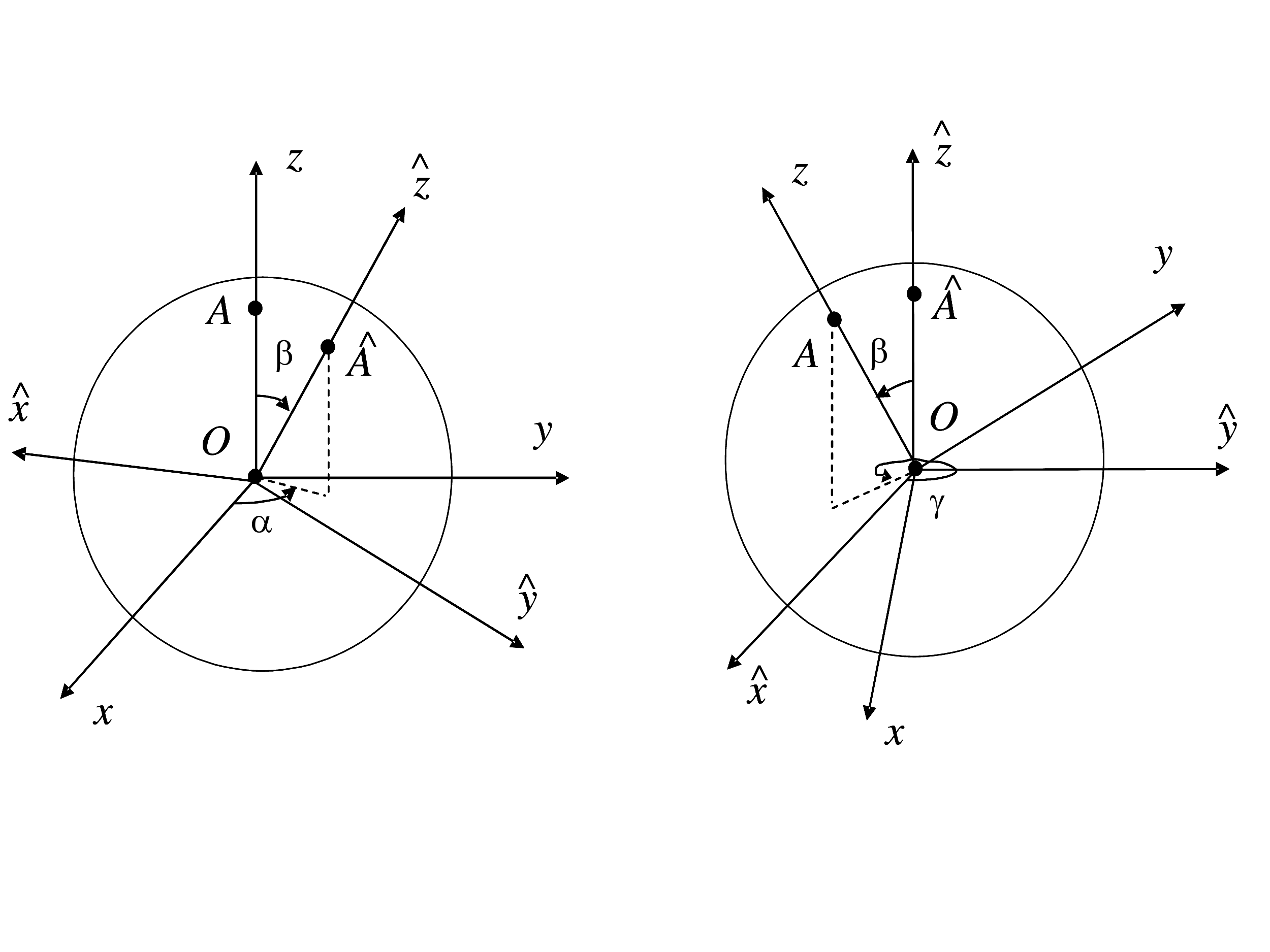}
\caption{Rotation angles $\protect\alpha ,\protect\beta ,$ and $\protect \gamma $ defined as spherical angles $(\protect\beta ,\protect\alpha)$ of the rotated $z$-axis in the original reference frame and
sperical angles $(\protect\beta ,\protect\gamma) $ of the
original $z$-axis in the rotated reference frame.}
\label{Fig1}
\end{figure}

An arbitrary rotation transform can be specified
by three Euler angles of rotation. We slightly modify these angles to be
consistent with the rotation angles $\alpha ,\beta
,\gamma $ defined in \cite{Gumerov05:Book}. Let $\mathbf{i}_{x},\mathbf{i}_{y},$ and $%
\mathbf{i}_{z}$ be the Cartesian basis vectors of the original reference
frame, while $\widehat{\mathbf{i}}_{x},\widehat{\mathbf{i}}_{y},$ and $%
\widehat{\mathbf{i}}_{z}$ be the respective basis vectors of the rotated
reference frame. Cartesian coordinates of $\widehat{\mathbf{i}}_{z}$ in the
original reference frame and $\mathbf{i}_{z}$ in the rotated reference frame
can be written as 
\begin{eqnarray}
\widehat{\mathbf{i}}_{z} &=&\left( \sin \beta \cos \alpha ,\sin \beta \sin
\alpha ,\cos \beta \right) ,  \label{r1} \\
\mathbf{i}_{z} &=&\left( \sin \beta \cos \gamma ,\sin \beta \sin \gamma
,\cos \beta \right) .  \notag
\end{eqnarray}
Figure \ref{Fig1} illustrates the rotation angles and reference frames. Note
that  
\begin{equation}
Q^{-1}\left( \alpha ,\beta ,\gamma \right) =Q^{T}\left( \alpha ,\beta
,\gamma \right) =Q\left( \gamma ,\beta ,\alpha \right) ,  \label{r2}
\end{equation}
where $Q$ is the rotation matrix and superscript $T$ denotes transposition.
We also note that as $\beta $ is related to spherical angle $\theta $, so
its range can be limited by half period $\beta \in \left[ 0,\pi \right] $,
while for $\alpha $ and $\gamma $ we have full periods $\alpha \in \left[
0,2\pi \right) ,$ $\gamma \in \left[ 0,2\pi \right) $. 

Let us introduce the Euler rotation angles, $\alpha _{E},\beta _{E},$ and $\gamma _{E}$, where general rotation is defined as rotation around original $z$ axis by angle $\alpha _{E}$ followed by rotation around about the new $y$ axis by angle $\beta _{E}$ and, finally, by rotation around the new $z$ axis by angle $\gamma _{E}$, then angles $\alpha ,\beta ,\gamma $ are simply related to that as 
\begin{equation}
\alpha =\alpha _{E},\quad \beta =\beta _{E},\quad \gamma =\pi -\gamma _{E}.
\label{r2.1}
\end{equation}%
Note that in \cite{Gumerov05:Book}, the Euler angles were introduced
differently ($\alpha =\pi -\alpha _{E},\quad \beta =\beta _{E},\quad \gamma
=\gamma _{E}$), so formulae obtained via such decomposition should be
modified if the present work is to be combined with those relations. Elementary rotation
matrices about axes $z$ and $y$ are 
\begin{equation}
Q_{z}\left( \alpha _{E}\right) =\left( 
\begin{array}{ccc}
\cos \alpha _{E} & \sin \alpha _{E} & 0 \\ 
-\sin \alpha _{E} & \cos \alpha _{E} & 0 \\ 
0 & 0 & 1%
\end{array}%
\right) ,\quad Q_{y}\left( \beta _{E}\right) =\left( 
\begin{array}{ccc}
\cos \beta _{E} & 0 & -\sin \beta _{E} \\ 
0 & 1 & 0 \\ 
\sin \beta _{E} & 0 & \cos \beta _{E}%
\end{array}%
\right) .  \label{r2.1.1}
\end{equation}%
The standard Euler rotation matrix decomposition $Q=Q_{z}\left( \gamma
_{E}\right) Q_{y}\left( \beta _{E}\right) Q_{z}\left( \alpha _{E}\right) $
turns to 
\begin{equation}
Q=Q_{z}\left( \pi -\gamma \right) Q_{y}\left( \beta \right) Q_{z}\left(
\alpha \right) .  \label{r2.1.2}
\end{equation}%
More symmetric forms with respect to angle $\beta $ can be obtained, if we
introduce elementary matrices $A$ and $B$ as follows 
\begin{eqnarray}
A\left( \gamma \right) &=&\left( 
\begin{array}{ccc}
\sin \gamma & \cos \gamma & 0 \\ 
-\cos \gamma & \sin \gamma & 0 \\ 
0 & 0 & 1%
\end{array}%
\right) =Q_{z}\left( \frac{\pi }{2}-\gamma \right) ,  \label{r2.1.3} \\
B\left( \beta \right) &=&\left( 
\begin{array}{ccc}
-1 & 0 & 0 \\ 
0 & -\cos \beta & \sin \beta \\ 
0 & \sin \beta & \cos \beta%
\end{array}%
\right) =Q_{z}\left( \frac{\pi }{2}\right) Q_{y}\left( \beta \right)
Q_{z}\left( \frac{\pi }{2}\right) ,  \notag
\end{eqnarray}%
which results in decomposition 
\begin{equation}
Q\left( \alpha ,\beta ,\gamma \right) =A\left( \gamma \right) B\left( \beta
\right) A^{T}\left( \alpha \right) .  \label{r2.2}
\end{equation}
The property of rotation transform is that any subspace of degree $n$ is
transformed independently of a subspace of a different degree. Also rotation
around the $z$ axis results in  diagonal  rotation transform operator, as
it is seen from the definition provided by Eq. (\ref{3}). So, the rotation
transform can be described as 
\begin{equation}
\widehat{C}_{n}^{m^{\prime }}=\sum_{m=-n}^{n}T_{n}^{m^{\prime }m}\left(
\alpha ,\beta ,\gamma \right) C_{n}^{m},\quad T_{n}^{m^{\prime
}m}=e^{-im^{\prime }\gamma }H_{n}^{m^{\prime }m}\left( \beta \right)
e^{im\alpha },  \label{r3}
\end{equation}%
where $\widehat{C}_{n}^{m^{\prime }}$ are the expansion coefficients of
function $f$ over the spherical harmonic basis in the rotated reference
frame, 
\begin{equation}
f\left( \theta ,\varphi \right)
=\sum_{n=0}^{p-1}\sum_{m=-n}^{n}C_{n}^{m}Y_{n}^{m}\left( \theta ,\varphi
\right) =\sum_{n=0}^{p-1}\sum_{m^{\prime }=-n}^{n}\widehat{C}_{n}^{m^{\prime
}}Y_{n}^{m^{\prime }}\left( \widehat{\theta },\widehat{\varphi }\right) =%
\widehat{f}\left( \widehat{\theta },\widehat{\varphi }\right) .  \label{r4}
\end{equation}%
Particularly, if $C_{n}^{m}=\delta _{\nu m},$ where $\delta _{\nu m}$ is
Kronecker's delta, we have from Eqs (\ref{r3}) and (\ref{r4}) for each
subspace 
\begin{equation}
Y_{n}^{\nu }\left( \theta ,\varphi \right) =e^{i\nu \alpha }\sum_{m^{\prime
}=-n}^{n}H_{n}^{m^{\prime }\nu }\left( \beta \right) e^{-im^{\prime }\gamma
}Y_{n}^{m^{\prime }}\left( \widehat{\theta },\widehat{\varphi }\right)
,\quad \nu =-n,...,n.  \label{r4.1}
\end{equation}

It should be mentioned then that the matrix with elements $T_{n}^{m^{\prime
}m}\left( \alpha ,\beta ,\gamma \right) $, which we denote as $\mathbf{Rot}%
\left( Q\left( \alpha ,\beta ,\gamma \right) \right) $ and its invariant
subspaces as\linebreak\ $\mathbf{Rot}_{n}(Q\left( \alpha ,\beta ,\gamma
\right) $, is the Wigner D-matrix in an irreducible representation of the
group of rigid body rotations SO(3) \cite{Wigner31:Book}, \cite%
{Vilenkin68:Book} (with slight modifications presented below). Particularly,
we have decompositions 
\begin{eqnarray}
\!\!\!\!\!\!\!\!\!\!\mathbf{Rot}_{n}\left( Q\left( \alpha ,\beta ,\gamma
\right) \right) \!\! &=&\!\!\mathbf{Rot}_{n}\left( Q_{z}\left( \pi -\gamma
\right) \right) \mathbf{Rot}_{n}\left( Q_{y}\left( \beta \right) \right) 
\mathbf{Rot}_{n}\left( Q_{z}\left( \alpha \right) \right)  \label{r4.1.1} \\
&=&\mathbf{Rot}_{n}\left( A\left( \gamma \right) \right) \mathbf{Rot}%
_{n}\left( B\left( \beta \right) \right) \mathbf{Rot}_{n}\left( A\left(
-\alpha \right) \right) .  \notag
\end{eqnarray}%
Since $Q_{y}\left( 0\right) $ and $Q_{z}\left( 0\right) $ are identity
matrices (see Eq. (\ref{r2.1.1})), then corresponding matrices $\mathbf{Rot}%
\left( Q_{y}\left( 0\right) \right) $ and $\mathbf{Rot}\left( Q_{y}\left(
0\right) \right) $ are also identity matrices. So, from Eqs (\ref{r4.1.1})
and (\ref{r3}) we obtain 
\begin{eqnarray}
\mathbf{Rot}_{n}\left( Q_{z}\left( \alpha \right) \right) &=&\mathbf{Rot}%
_{n}\left( Q\left( \alpha ,0,\pi \right) \right) =\left\{ \left( -1\right)
^{m^{\prime }}H_{n}^{m^{\prime }m}\left( 0\right) e^{im\alpha }\right\} ,
\label{r4.1.2} \\
\mathbf{Rot}_{n}\left( Q_{y}\left( \beta \right) \right) &=&\mathbf{Rot}%
_{n}\left( Q\left( 0,\beta ,\pi \right) \right) =\left\{ \left( -1\right)
^{m^{\prime }}H_{n}^{m^{\prime }m}\left( \beta \right) \right\} ,  \notag
\end{eqnarray}%
where in the figure brackets we show the respective elements of the
matrices. Since $\mathbf{Rot}_{n}\left( Q_{y}\left( 0\right) \right) $ is
the identity matrix, the latter equation provides 
\begin{equation}
H_{n}^{m^{\prime }m}\left( 0\right) =\left( -1\right) ^{m^{\prime }}\delta
_{m^{\prime }m}.  \label{r4.1.3}
\end{equation}%
Being used in the first equation (\ref{r4.1.2}) this results in 
\begin{equation}
\mathbf{Rot}_{n}\left( Q_{z}\left( \alpha \right) \right) =\mathbf{Rot}%
_{n}\left( Q\left( \alpha ,0,\pi \right) \right) =\left\{ e^{im\alpha
}\delta _{m^{\prime }m}\right\} ,  \label{r4.1.4}
\end{equation}%
which shows that 
\begin{equation}
\mathbf{Rot}_{n}\left( A\left( \gamma \right) \right) =\mathbf{Rot}%
_{n}\left( Q_{z}\left( \frac{\pi }{2}-\gamma \right) \right) =\left\{
e^{im\pi /2}e^{-im\gamma }\delta _{m^{\prime }m}\right\} .  \label{r4.1.5}
\end{equation}%
We also have immediately from Eq. (\ref{r3}) 
\begin{equation}
\mathbf{Rot}_{n}\left( B\left( \beta \right) \right) =\mathbf{Rot}_{n}\left(
Q\left( 0,\beta ,0\right) \right) =\left\{ H_{n}^{m^{\prime }m}\left( \beta
\right) \right\} .  \label{r4.1.6}
\end{equation}

The rotation coefficients $H_{n}^{m^{\prime }m}\left( \beta \right) $ are
real and are simply related to the Wigner's (small) d-matrix elements ($%
d_{n}^{m^{\prime }m}\left( \beta \right) ,$ elsewhere, e.g. \cite%
{Stein61:QAM}),

\begin{equation}
d_{n}^{m^{\prime }m}\left( \beta \right) =\left( -1\right) ^{m^{\prime
}-m}\rho _{n}^{m^{\prime
}m}\!\!\!\!\!\!\!\!\!\!\!\!\!\!\!\!\!\!\!\!\sum_{\sigma =\max \left(
0,-(m^{\prime }-m)\right) }^{\min \left( n-m^{\prime },n+m\right) }\frac{%
\left( -1\right) ^{\sigma }\cos ^{2n-2\sigma +m-m^{\prime }}\frac{1}{2}\beta
\sin ^{2\sigma +m^{\prime }-m}\frac{1}{2}\beta }{\sigma !\left( n+m-\sigma
\right) !\left( n-m^{\prime }-\sigma \right) !\left( m^{\prime }-m+\sigma
\right) !},  \label{r4.2}
\end{equation}%
but slightly different due to the difference in definition of spherical
harmonics and rotation matrix. In this expression we defined $\rho
_{n}^{m^{\prime }m}$ as%
\begin{equation}
\rho _{n}^{m^{\prime }m}=\left[ (n+m)!(n-m)!(n+m^{\prime })!(n-m^{\prime })!%
\right] ^{1/2}.  \label{r4.2.1}
\end{equation}%
Explicit expression for coefficients $H_{n}^{m^{\prime }m}\left( \beta
\right) $ can be obtained from Wigner's formula,

\begin{equation}
H_{n}^{m^{\prime }m}\left( \beta \right) =\epsilon _{m^{\prime }}\epsilon
_{m}\rho _{n}^{m^{\prime }m}\sum_{\sigma =\max \left( 0,-(m^{\prime
}+m)\right) }^{\min \left( n-m^{\prime },n-m\right) }\left( -1\right)
^{n-\sigma }h_{n}^{m^{\prime }m\sigma }\left( \beta \right) ,  \label{r5}
\end{equation}%
where 
\begin{equation}
h_{n}^{m^{\prime }m\sigma }\left( \beta \right) =\frac{\cos ^{2\sigma
+m+m^{\prime }}\frac{1}{2}\beta \sin ^{2n-2\sigma -m-m^{\prime }}\frac{1}{2}%
\beta }{\sigma !\left( n-m^{\prime }-\sigma \right) !\left( n-m-\sigma
\right) !\left( m^{\prime }+m+\sigma \right) !},  \label{r5.1}
\end{equation}%
and symbol $\epsilon _{m^{\prime }}$ is defined by Eq. (\ref{r6}). Note that
summation limits in Eq. (\ref{r5}) can look a bit complicated, but this can
be avoided, if we simply define $1/(-n)!=0$ for $n=1,2,...$ (which is
consistent with the limit of $1/\Gamma (-n)$ for $n=0,1,...$, where $\Gamma $
is the gamma-function), so 
\begin{equation}
H_{n}^{m^{\prime }m}\left( \beta \right) =\epsilon _{m^{\prime }}\epsilon
_{m}\rho _{n}^{m^{\prime }m}\sum_{\sigma =-\infty }^{\infty }\left(
-1\right) ^{n-\sigma }h_{n}^{m^{\prime }m\sigma }\left( \beta \right) .
\label{r7}
\end{equation}%
The Wigner's (small) d-matrix elements are related to $H_{n}^{m^{\prime
}m}\left( \beta \right) $ coefficients as 
\begin{equation}
d_{n}^{m^{\prime }m}\left( \beta \right) =\epsilon _{m^{\prime }}\epsilon
_{-m}H_{n}^{m^{\prime }m}\left( \beta \right) ,  \label{r8}
\end{equation}%
which can be checked directly using Eqs (\ref{r4.2}) and (\ref{r5}) and
symmetries (\ref{s1}).

\subsection{Symmetries}

There are several symmetries of the rotation coefficients, from which the
following are important for the proposed algorithm 
\begin{eqnarray}
H_{n}^{m^{\prime }m}\left( \beta \right) &=&H_{n}^{mm^{\prime }}\left( \beta
\right) ,\quad  \label{s1} \\
H_{n}^{m^{\prime }m}\left( \beta \right) &=&H_{n}^{-m^{\prime },-m}\left(
\beta \right) ,  \notag \\
H_{n}^{m^{\prime }m}\left( \pi -\beta \right) &=&(-1)^{n+m^{\prime
}+m}H_{n}^{-m^{\prime }m}\left( \beta \right) ,  \notag \\
H_{n}^{m^{\prime }m}\left( -\beta \right) &=&(-1)^{m^{\prime
}+m}H_{n}^{m^{\prime }m}\left( \beta \right) .  \notag
\end{eqnarray}

The first symmetry follows trivially from Eq. (\ref{r5}), which is symmetric
with respect to $m^{\prime }$ and $m$.

The second symmetry also can be proved using Eq. (\ref{r5}) or its analog (%
\ref{r7}). This can be checked straightforward using $\epsilon _{-m^{\prime
}}=(-1)^{m^{\prime }}\epsilon _{m^{\prime }}$ and replacement $\sigma
=\sigma ^{\prime }-m^{\prime }-m$ in the sum. The third symmetry (\ref{s1})
can be also obtained from Eq. (\ref{r7}) using $\sin \frac{1}{2}\left( \pi
-\beta \right) =\cos \frac{1}{2}\beta $, $\epsilon _{-m^{\prime
}}=(-1)^{m^{\prime }}\epsilon _{m^{\prime }}$ and replacement $\sigma
=n-\sigma ^{\prime }-m$ in the sum. The fourth symmetry follows simply from
Eq. (\ref{r7}). It is not needed for $\beta \in \left[ 0,\pi \right] $, but
we list it here for completness, as full period change of $\beta $ sometimes
may be needed.

\subsection{Particular values}

For some values of $m^{\prime },m,$ and $n$ coefficients $H_{n}^{m^{\prime
}m}$ can be computed with minimal cost, which does not require summation (%
\ref{r5}). For example, the addition theorem for spherical harmonics can be
written in the form 
\begin{equation}
P_{n}\left( \cos \theta \right) =\frac{4\pi }{2n+1}\sum_{m^{\prime
}=-n}^{n}Y_{n}^{-m^{\prime }}\left( \theta _{2},\varphi _{2}\right)
Y_{n}^{m^{\prime }}\left( \theta _{1},\varphi _{1}\right) ,  \label{v1}
\end{equation}%
where $\theta $ is the angle between points on a unit sphere with spherical
coordinates $\left( \theta _{1},\varphi _{1}\right) $ and $\left( \theta
_{2},\varphi _{2}\right) $. From definition of rotation angles $\alpha
,\beta ,\gamma $, we can see then that if $\left( \theta _{1},\varphi
_{1}\right) =\left( \widehat{\theta },\widehat{\varphi }\right) $ are
coordinates of the point in the rotated reference frame, which coordinates
in the original reference frame are $\left( \theta ,\varphi \right) $ then $%
\left( \theta _{2},\varphi _{2}\right) =\left( \beta ,\gamma \right) $,
since the scalar product of radius-vectors pointed to $\left( \widehat{%
\theta },\widehat{\varphi }\right) $ and $\left( \beta ,\gamma \right) $
(the $\mathbf{i}_{z}$ in the rotated frame, Eq. (\ref{r1})) will be $\cos
\theta $ (the $z$-coordinate of the same point in the original reference
frame). Comparing this with Eq. (\ref{r4.1}) for $\nu =0$ and using
definition of spherical harmonics (\ref{3}) and (\ref{4}), we obtain 
\begin{equation}
H_{n}^{m^{\prime }0}\left( \beta \right) =(-1)^{m^{\prime }}\sqrt{\frac{%
(n-\left\vert m^{\prime }\right\vert )!}{(n+\left\vert m^{\prime
}\right\vert )!}}P_{n}^{\left\vert m^{\prime }\right\vert }(\cos \beta ),
\label{v2}
\end{equation}%
since the relation is valid for arbitrary point on the sphere. Note that
computation of normalized associated Legendre function can be done using
well-known stable recursions and standard library routines are available.
This results in $O\left( p^{2}\right) $ cost of computation of all $%
H_{n}^{m^{\prime }0}\left( \beta \right) $ for $n=0,...,p-1,$ $m^{\prime
}=-n,...,n$.

Another set of easily and accurately computable values comes from Wigner's
formula (\ref{r5}), where the sum reduces to a single term ($\sigma =0)$ at $%
n=m$. We have in this case 
\begin{equation}
H_{n}^{m^{\prime }n}\left( \beta \right) =\epsilon _{m^{\prime }}\left[ 
\frac{(2n)!}{\left( n-m^{\prime }\right) !\left( n+m^{\prime }\right) !}%
\right] ^{1/2}\cos ^{n+m^{\prime }}\frac{1}{2}\beta \sin ^{n-m^{\prime }}%
\frac{1}{2}\beta .  \label{v3}
\end{equation}

\subsection{Axis flip transform}

There exist more expressions for $H_{n}^{m^{\prime }m}\left( \beta \right) $
via finite sums and the present algorithm uses one of those. This relation
can be obtained from consideration of the axis flip transform. A composition
of rotations which puts axis $y$ in position of axis $z$ and then performs
rotation about the $z$ axis, which is described by diagonal matrix followed
by the inverse transform is well-known and used in some algorithms. The flip
transform can be described by the following formula 
\begin{equation}
Q_{y}\left( \beta \right) =Q_{z}\left( -\frac{\pi }{2}\right) Q_{y}\left( -%
\frac{\pi }{2}\right) Q_{z}\left( \beta \right) Q_{y}\left( \frac{\pi }{2}%
\right) Q_{z}\left( \frac{\pi }{2}\right) ,  \label{af1}
\end{equation}
which meaning in rather obvious from geometry as it decomposes rotation
around the $y$-axis into $\pi /2$ rotations about $z$ and $y$ axes, rotation
by angle $\beta $ around $z$-axis followed by the inverse rotations about $y$
and $z$ axes.

From Eq. (\ref{r2.1.3}) we have 
\begin{eqnarray}
Q_{y}\left( \frac{\pi }{2}\right) &=&Q_{z}\left( -\frac{\pi }{2}\right)
B\left( \frac{\pi }{2}\right) Q_{z}\left( -\frac{\pi }{2}\right) ,
\label{af2} \\
Q_{y}\left( -\frac{\pi }{2}\right) &=&Q_{y}^{T}\left( \frac{\pi }{2}\right)
=Q_{z}\left( \frac{\pi }{2}\right) B\left( \frac{\pi }{2}\right) Q_{z}\left( 
\frac{\pi }{2}\right) .  \notag
\end{eqnarray}
Using the same equation one can express $Q_{y}\left( \beta \right) $ via $%
B\left( \beta \right) $ and obtain from Eqs (\ref{af1}) and (\ref{af2})

\begin{equation}
B\left( \beta \right) =Q_{z}\left( \frac{\pi }{2}\right) B\left( \frac{\pi }{%
2}\right) Q_{z}\left( \beta \right) B\left( \frac{\pi }{2}\right)
Q_{z}\left( \frac{\pi }{2}\right) .  \label{af3}
\end{equation}%
The representation of this decomposition for each subspace $n$ results in 
\begin{eqnarray}
\mathbf{Rot}_{n}\left( B\left( \beta \right) \right) &=&\mathbf{Rot}%
_{n}\left( Q_{z}\left( \frac{\pi }{2}\right) \right) \mathbf{Rot}_{n}\left(
B\left( \frac{\pi }{2}\right) \right) \times  \label{af4} \\
&&\mathbf{Rot}_{n}\left( Q_{z}\left( \beta \right) \right) \mathbf{Rot}%
_{n}\left( B\left( \frac{\pi }{2}\right) \right) \mathbf{Rot}_{n}\left(
Q_{z}\left( \frac{\pi }{2}\right) \right) .  \notag
\end{eqnarray}%
Using expressions (\ref{r4.1.4}) and (\ref{r4.1.6}), we can rewrite this
relation in terms of matrix elements 
\begin{equation}
H_{n}^{m^{\prime }m}\left( \beta \right) =\sum_{\nu =-n}^{n}e^{im^{\prime
}\pi /2}H_{n}^{m^{\prime }\nu }\left( \frac{\pi }{2}\right) e^{i\nu \beta
}H_{n}^{\nu m}\left( \frac{\pi }{2}\right) e^{im\pi /2}.  \label{af5}
\end{equation}%
Since $H_{n}^{m^{\prime }m}\left( \beta \right) $ is real, we can take the
real part of the right hand side of this relation, to obtain 
\begin{equation}
H_{n}^{m^{\prime }m}\left( \beta \right) =\sum_{\nu =-n}^{n}H_{n}^{m^{\prime
}\nu }\left( \frac{\pi }{2}\right) H_{n}^{m\nu }\left( \frac{\pi }{2}\right)
\cos \left( \nu \beta +\frac{\pi }{2}\left( m^{\prime }+m\right) \right) ,
\label{af6}
\end{equation}%
where we used the first symmetry (\ref{s1}). Note also that the third
symmetry (\ref{s1}) applied to $H_{n}^{m^{\prime }m}\left( \pi /2\right) $
results in 
\begin{eqnarray}
H_{n}^{m^{\prime }m}\left( \beta \right) &=&H_{n}^{m^{\prime }0}\left( \frac{%
\pi }{2}\right) H_{n}^{m0}\left( \frac{\pi }{2}\right) \cos \left( \frac{\pi 
}{2}\left( m^{\prime }+m\right) \right) +  \label{af7} \\
&&2\sum_{\nu =1}^{n}H_{n}^{m^{\prime }\nu }\left( \frac{\pi }{2}\right)
H_{n}^{m\nu }\left( \frac{\pi }{2}\right) \cos \left( \nu \beta +\frac{\pi }{%
2}\left( m^{\prime }+m\right) \right) .  \notag
\end{eqnarray}

\subsection{Recursions}
Several recursions for computation of
coefficients $H_{n}^{m^{\prime }m}\left( \beta \right) $ were derived from the
invariancy of differential operator $\nabla $ \cite{Gumerov03:SISC}, including the following one, which was suggested for
computing all $H_{n}^{m^{\prime }m}\left( \beta \right) $%
\begin{eqnarray}
b_{n}^{m}H_{n-1}^{m^{\prime },m+1} &=&\frac{1-\cos \beta }{2}%
b_{n}^{-m^{\prime }-1}H_{n}^{m^{\prime }+1,m}-  \label{rec1} \\
&&\frac{1+\cos \beta }{2}b_{n}^{m^{\prime }-1}H_{n}^{m^{\prime }-1,m}-\sin
\beta a_{n-1}^{m^{\prime }}H_{n}^{m^{\prime }m},  \notag
\end{eqnarray}%
where $n=2,3,...,\quad m^{\prime }=-n+1,...,n-1,\quad m=0,...,n-2,$ and
\begin{eqnarray}
a_{n}^{m} &=&a_{n}^{-m}=\sqrt{\frac{(n+1+m)(n+1-m)}{\left( 2n+1\right)
\left( 2n+3\right) }},\quad n\geqslant \left\vert m\right\vert .\qquad 
\label{rec1.1} \\
a_{n}^{m} &=&0,\quad n<\left\vert m\right\vert ,  \notag
\end{eqnarray}%
\begin{equation}
\!\!\!\!\!\!\!\!\!\!\!\!\!b_{n}^{m}\! =\text{sgn}\left( m\right)
\!\sqrt{\frac{(n-m-1)(n-m)}{\left( 2n-1\right) \left( 2n+1\right) }},\quad
n\geqslant \left\vert m\right\vert \quad \!\!\!\!
\quad \!n<\left\vert m\right\vert \!,  
\label{rec1.2} 
\end{equation}
\begin{equation}
\text{sgn}\left( m\right) =\left\{ 
\begin{array}{c}
1,\quad m\geqslant 0 \\ 
-1,\quad m<0%
\end{array}%
\right. .  \label{rec1.3}
\end{equation}%
This recursion allows one to get  $\left\{ H_{n}^{m^{\prime
},m+1}\right\} $ from $\left\{ H_{n}^{m^{\prime }m}\right\}$.
Once the value for $m=0,$ $H_{n}^{m^{\prime }0}$, is known from Eq. (\ref{v2}%
) for any $n$ and $m^{\prime }$ all coefficients can be computed. The
cost of the procedure is $O\left( p^{3}\right) $ as $n$ is
limited by $p-1$ (to get that $H_{n}^{m^{\prime }0}$ should be computed up
to $n=2p-2$). The relation was extensively tested and used in
the FMM (e.g. \cite{Gumerov05:JASA}), however, it showed
numerical instability $p\gtrsim 100.$

From the commutativity of rotations around the axis $y$,
\begin{equation}
Q_{y}\left( \beta \right) Q_{y}^{\prime }\left( 0\right) =Q_{y}^{\prime
}\left( 0\right) Q_{y}\left( \beta \right), \quad Q_{y}^{\prime }\left( \beta \right) =\frac{dQ_{y}}{d\beta }.  \label{rec1.4}
\end{equation}
one may derive a recurrence relation \cite{Gumerov05:Book}. For the $n$th invariant rotation transform subspace, the relation (\ref {rec1.4}) gives 
\begin{equation}
\mathbf{Rot}_{n}\left( Q_{y}\left( \beta \right) \right) \mathbf{Rot}%
_{n}\left( Q_{y}^{\prime }\left( 0\right) \right) =\mathbf{Rot}_{n}\left(
Q_{y}^{\prime }\left( 0\right) \right) \mathbf{Rot}_{n}\left( Q_{y}\left(
\beta \right) \right) .  \label{rec1.6}
\end{equation}%
Differentiating the r.h.s. of (\ref{r7}) w.r.to $\beta$ and evaluating at $\beta =0$ we get 
\begin{equation}
\left. \frac{dH_{n}^{m^{\prime }m}\left( \beta \right) }{d\beta }\right\vert
_{\beta =0}=c_{n}^{m^{\prime }-1}\delta _{m,m^{\prime }-1}+c_{n}^{m^{\prime
}}\delta _{m,m^{\prime }+1},  \label{rec1.7}
\end{equation}%
\begin{equation}
\!\!\!\!\!\!\!\!\!\!\!\!\!\!\!\!\!\!c_{n}^{m}\!=\!\frac{1}{2}(-1)^{m}\text{sgn}\left( m\right) \left[ (n-m)(n+m+1)\right]^{1/2},\!\quad\! m=-n-1,...,n.\!  \label{rec1.8}
\end{equation}%
Using (\ref{r4.1.2}),  relation (\ref{rec1.6}) can be
rewritten as 
\begin{equation}
\!\!\!\!\!\!\!\!\!\!\!\!\!\sum_{\nu =-n}^{n}H_{n}^{m^{\prime }\nu }\left( \beta \right)
\left( -1\right) ^{\nu }\left. \frac{dH_{n}^{\nu m}\left( \beta \right) }{%
d\beta }\right\vert _{\beta =0}\!\!\!\!\!=\!\!\!\!\sum_{\nu =-n}^{n}\left. 
\frac{dH_{n}^{m^{\prime }\nu }\left( \beta \right) }{d\beta }\right\vert
_{\beta =0}\!\!\!\!\!\!\!\!\!\!\left( -1\right) ^{\nu }H_{n}^{\nu m}\left( \beta \right)
.  \label{rec1.9}
\end{equation}%
Using Eq. (\ref{rec1.7}) we obtain the recurrence relation for $H_{n}^{m^{\prime },m}( \beta)$
\begin{equation}
\!\!\!\!\!\!\!\!\!\!\!\!d_{n}^{m-1}H_{n}^{m^{\prime },m-1}
-d_{n}^{m}H_{n}^{m^{\prime },m+1} =d_{n}^{m^{\prime}-1}H_{n}^{m^{\prime }-1,m} -d_{n}^{m^{\prime }}H_{n}^{m^{\prime }+1,m},  \label{rec1.10}
\end{equation}%
\begin{equation}
\!\!\!\!\!\!\!\!\!\!\!\!\!\!\!\!\!d_{n}^{m}=(-1)^{m}c_{n}^{m}=\frac{\text{sgn}(m)}{2} \left[
(n-m)(n+m+1)\right] ^{1/2},\quad m=-n-1,...,n.  \label{rec1.11}
\end{equation}
 {\bf This recurrence is one of the main results of this paper}. In contrast to the recurrence (\ref{rec1}), (\ref{rec1.10}) relates
values of rotation coefficients\ $H_{n}^{m^{\prime }m}$ {\em within the same
subspace} $n$. Thus, if boundary values for a subspace are provided, all other coefficients can be found.
\section{Bounds for rotation coefficients}
It should be noticed that the $\left( 2n+1\right) \times $ $\left( 2n+1\right) $
matix $\mathbf{H}_{n}\left( \beta \right) =\mathbf{Rot}_{n}\left( B\left(
\beta \right) \right) $ with entries $H_{n}^{m^{\prime }m}\left( \beta
\right) $, $m^{\prime },m=-n,...,n$ is real, unitary, and self-adjoint
(Hermitian). This follows from 
\begin{equation}
\left[ \mathbf{H}_{n}\left( \beta \right) \right] ^{2}=\mathbf{I}_{n}\mathbf{%
,\quad H}_{n}\left( \beta \right) =\mathbf{H}_{n}^{T}\left( \beta \right) ,
\label{b1}
\end{equation}%
where $\mathbf{I}_{n}$ is $\left( 2n+1\right) \times $ $\left( 2n+1\right) $
identity matrix. Particularly, this shows that the norm of matrix $\mathbf{H}%
_{n}\left( \beta \right) $ is unity and 
\begin{equation}
\left\vert H_{n}^{m^{\prime }m}\left( \beta \right) \right\vert \leqslant
1,\quad n=0,1,...,\quad m^{\prime },m=-n,...,n.  \label{b2}
\end{equation}

While bound (\ref{b2}) is applicable for any values of $\beta ,n,m$ and $%
m^{\prime }$, we note that in certain regions of the parameter space it can
be improved.

\subsection{General bound}

To get such a bound we note that due to the third symmetry (\ref{s1}) it is
sufficient to consider $\beta $ in range $0\leqslant \beta \leqslant \pi /2$%
. The first and the second symmetries (\ref{s1}) provide that only
non-negative $m$ can be considered, $m\geqslant 0$, and also $m^{\prime }$
from the range $\left\vert m^{\prime }\right\vert \leqslant m$. The latter
provides $m+m^{\prime }\geqslant 0,$ $m^{\prime }\leqslant m,$ and $%
n-m^{\prime }\geqslant n-m.$ Hence in this range Eq. (\ref{r5}) can be
written in the form 
\begin{equation}
H_{n}^{m^{\prime }m}\left( \beta \right) =\epsilon _{m^{\prime }}\epsilon
_{m}\rho _{n}^{m^{\prime }m}\sum_{\sigma =0}^{n-m}\left( -1\right)
^{n-\sigma }h_{n}^{m^{\prime }m\sigma }\left( \beta \right) ,  \label{b3}
\end{equation}%
where $h_{n}^{m^{\prime }m\sigma }\left( \beta \right) \geqslant 0$, and we
can bound $\left\vert H_{n}^{m^{\prime }m}\left( \beta \right) \right\vert $
as 
\begin{eqnarray}
\left\vert H_{n}^{m^{\prime }m}\left( \beta \right) \right\vert &\leqslant
&\rho _{n}^{m^{\prime }m}\sum_{\sigma =0}^{n-m}h_{n}^{m^{\prime }m\sigma
}\left( \beta \right)  \label{b4} \\
&\leqslant &\rho _{n}^{m^{\prime }m}(n-m+1)h_{n}^{m^{\prime }ms}\left( \beta
\right) ,\text{ }  \notag
\end{eqnarray}%
where 
\begin{equation}
h_{n}^{m^{\prime }ms}\left( \beta \right) =\max_{0\leqslant \sigma \leqslant
n-m}h_{n}^{m^{\prime }m\sigma }\left( \beta \right) .  \label{b5}
\end{equation}%
In other words $s$ is the value of $\sigma $ at which $h_{n}^{m^{\prime
}m\sigma }\left( \beta \right) $ achieves its maximum at given $\beta ,n,m$
and $m^{\prime }$.

Note then that for $n=m$ we have only one term, $\sigma =0,$ the sum (\ref{b3}%
), so $s=0$ and in this case $\left| H_{n}^{m^{\prime }m}\left( \beta
\right) \right| $ reaches the bound (\ref{b4}), as it also follows from Eqs (%
\ref{r5.1}) and (\ref{v3}). Another simple case is realized for $\beta =0.$
In this case, again, the sum has only one non-zero term at $\sigma =n-m$ and 
\begin{equation}
h_{n}^{m^{\prime }m,n-m}\left( 0\right) =\frac{\delta _{m^{\prime }m}}{%
(n-m)!(n+m)!}.  \label{b6}
\end{equation}
While exact value in this case is $\left| H_{n}^{m^{\prime }m}\left(
0\right) \right| =\delta _{m^{\prime }m}$ (see also Eq. (\ref{r4.1.3})),
bound (\ref{b4}) at $s=n-m$ provides $\left| H_{n}^{m^{\prime }m}\left(
\beta \right) \right| \leqslant \left( n-m+1\right) \delta _{m^{\prime }m}$,
which is correct, but is not tight.

To find the maximum of $h_{n}^{m^{\prime }m\sigma }\left( \beta \right) $ for
$\beta \neq 0$ and $n-m>0$ (so, also $n-m^{\prime }>0$) we can consider the
ratio of the consequent terms in the sum 
\begin{equation}
r_{n}^{m^{\prime }m\sigma }\left( \beta \right) =\frac{h_{n}^{m^{\prime
}m,\sigma +1}\left( \beta \right) }{h_{n}^{m^{\prime }m\sigma }\left( \beta
\right) }=\frac{\left( n-m^{\prime }-\sigma \right) \left( n-m-\sigma
\right) }{t^{2}\left( \sigma +1\right) \left( m^{\prime }+m+\sigma +1\right) 
},\quad t=\tan \frac{\beta }{2},  \label{b7}
\end{equation}%
where parameter $t$ is varying in the range $0<t\leqslant 1$, as we consider 
$0<\beta \leqslant \pi /2.$ This ratio considered as a function of $\sigma $
monotonously decay from its value at $\sigma =0$ to zero at $\sigma =n-m$
(as the numerator is a decaying function, while the denominator is a growing
function ($n-m^{\prime }\geqslant n-m$). So, if $r_{n}^{m^{\prime }m0}\left(
\beta \right) \leqslant 1$ then the maximum of $h_{n}^{m^{\prime }m\sigma
}\left( \beta \right) $ is reached at $\sigma =0$, otherwise it can be found
from simultaneous equations $r_{n}^{m^{\prime }ms}\left( \beta \right) >1,$ $%
r_{n}^{m^{\prime }m,s+1}\left( \beta \right) \leqslant 1.$ To treat both
cases, we consider roots of equation $r_{n}^{m^{\prime }m\sigma }\left(
\beta \right) =1$, which turns to a quadratic equation with respect to $%
\sigma ,$%
\begin{equation}
\left( n-m^{\prime }-\sigma \right) \left( n-m-\sigma \right) =t^{2}\left(
\sigma +1\right) \left( m^{\prime }+m+\sigma +1\right) .  \label{b8}
\end{equation}%
If there is no real non-negative roots in range $0\leqslant \sigma \leqslant
n-m$ then $s=0$, otherwise $s$ should be the integer part of the smallest
root of Eq. (\ref{b8}), as there may exist only one root of equation $%
r_{n}^{m^{\prime }m\sigma }\left( \beta \right) =1$ in range $0\leqslant
\sigma \leqslant n-m$ (the largest root in this case is at $\sigma
>n-m^{\prime }\geqslant n-m$). So, we have 
\begin{eqnarray}
\sigma _{1} &=&\frac{2n+2t^{2}-\left( 1-t^{2}\right) \left( m+m^{\prime
}\right) -\sqrt{D}}{2\left( 1-t^{2}\right) },  \label{b9} \\
D &=&\left( m-m^{\prime }\right) ^{2}+t^{2}\left( 2\left(
2n^{2}-m^{2}-(m^{\prime })^{2}\right) +t^{2}(m+m^{\prime })^{2}+4\left(
2n+1\right) \right) .  \notag
\end{eqnarray}%
This shows that $D\geqslant 0$, so the root is anyway real. Note also that $%
t=1$ $\left( \beta =\pi /2\right) $ is a special case, since at this value
equation (\ref{b8}) degenerates to a linear equation, which has 
root 
\begin{equation}
\sigma _{1}=\frac{\left( n-m^{\prime }\right) \left( n-m\right) -\left(
m+m^{\prime }+1\right) }{2(n+1)},\quad \left( t=1\right) .  \label{b10}
\end{equation}%
Summarizing, we obtain the following expression for $s$%
\begin{equation}
s=\left\{ 
\begin{array}{c}
\left[ \sigma _{1}\right] ,\quad \sigma _{1}\geqslant 0 \\ 
0,\quad \sigma _{1}<0,%
\end{array}%
\right.   \label{b11}
\end{equation}%
where $\left[ {}\right] $ denotes the integer part, $\sigma _{1}\left(
n,m,m^{\prime },t\right) $ for $t\neq 1$ is provided by Eq. (\ref{b9}), and
its limiting value at $t=1$ is given by Eq. (\ref{b10}). Equations (\ref{b4}%
) and (\ref{r5.1}) then yield 
\begin{equation}
\left\vert H_{n}^{m^{\prime }m}\left( \beta \right) \right\vert \leqslant 
\frac{\rho _{n}^{m^{\prime }m}(n-m+1)\cos ^{2s+m+m^{\prime }}\frac{\beta }{2}%
\sin ^{2n-2s-m-m^{\prime }}\frac{\beta }{2}}{s!\left( n-m^{\prime }-s\right)
!\left( n-m-s\right) !\left( m+m^{\prime }+s\right) !}.\text{ }  \label{b12}
\end{equation}

\subsection{Asymptotic behavior}

The bounds for the magnitude of the rotation coefficients are important for
study of their behavior at large $n$, as at large enough $n$ recursions
demonstrate instabilities, while direct computations using the sums becomes
difficult due to factorials of large numbers. So, we are going to obtain
asymptotics of expression (\ref{b12}) for $\left| H_{n}^{m^{\prime }m}\left(
\beta \right) \right| $ at $n\rightarrow \infty $. We note then that for
this purpose, we consider scaling of parameters $m,m^{\prime },$ and $s$,
i.e. we introduce new variables 
\begin{equation}
\mu =\frac{m}{n},\quad \mu ^{\prime }=\frac{m^{\prime }}{n},\quad \xi =\frac{%
s}{n},  \label{ab1}
\end{equation}
which, as follows from the above consideration, are in the range $0\leqslant
\mu \leqslant 1,$ $-\mu \leqslant \mu ^{\prime }\leqslant \mu ,$ $0\leqslant
\xi \leqslant 1-\mu $. The asymptotics can be constructed assuming that
these parameters are fixed, while $n\rightarrow \infty .$

Note now, that Eq. (\ref{b12}) can be written in the form 
\begin{eqnarray}
\left\vert H_{n}^{m^{\prime }m}\left( \beta \right) \right\vert  &\leqslant
&(n-m+1)\left[ C_{n-m}^{s}C_{n-m^{\prime }}^{s}C_{n+m}^{m+m^{\prime
}+s}C_{n+m^{\prime }}^{m+m^{\prime }+s}\right] ^{1/2}\times   \label{ab3} \\
&&\cos ^{2s+m+m^{\prime }}\frac{1}{2}\beta \sin ^{2n-2s-m-m^{\prime }}\frac{1%
}{2}\beta ,  \notag
\end{eqnarray}%
where 
\begin{equation}
C_{q}^{l}=\frac{q!}{l!\left( q-l\right) !},  \label{ab4}
\end{equation}%
are the binomial coefficients. Consider asymptotics of $C_{an}^{bn}$, where $%
a$ and $b$ are fixed, $0<b<a$, and $n\rightarrow \infty $. Using the
inequality, valid for $x>0,$ (e.g. see \cite{Abramowitz1972}) 
\begin{equation}
\sqrt{2\pi }\exp \left( \frac{2x+1}{2}\ln x-x\right) <x!<\sqrt{2\pi }\exp
\left( \frac{2x+1}{2}\ln x-x+\frac{1}{12x}\right) .  \label{ab5}
\end{equation}%
We can find that 
\begin{equation}
C_{an}^{bn}<C_{ab}n^{-1/2}e^{\lambda _{ab}n},  \label{ab6}
\end{equation}%
where 
\begin{equation}
\lambda _{ab}=a\ln a-b\ln b-(a-b)\ln \left( a-b\right) ,  \label{ab7}
\end{equation}%
and the constant $C_{ab}$ is bounded as 
\begin{equation}
C_{ab}\leqslant \sqrt{\frac{a}{2\pi b(a-b)}}e^{1/12}.  \label{ab8}
\end{equation}%
Using this bound in Eq. (\ref{ab3}) and definition (\ref{ab1}), we obtain 
\begin{equation}
\left\vert H_{n}^{m^{\prime }m}\left( \beta \right) \right\vert \leqslant
Ce^{\lambda n},  \label{ab9}
\end{equation}%
where $C$ is some constant depending on $\mu ,\mu ^{\prime },$and $\xi $,
while for $\lambda $ we have the following expression. 
\begin{eqnarray}
\lambda  &=&\frac{1-\mu }{2}\ln \left( 1-\mu \right) +\frac{1-\mu ^{\prime }%
}{2}\ln \left( 1-\mu ^{\prime }\right) +\frac{1+\mu }{2}\ln \left( 1+\mu
\right) +  \label{ab10} \\
&&\frac{1+\mu ^{\prime }}{2}\ln \left( 1+\mu ^{\prime }\right) -\xi \ln \xi
-\left( \mu +\mu ^{\prime }+\xi \right) \ln \left( \mu +\mu ^{\prime }+\xi
\right) -  \notag \\
&&\left( 1-\mu -\xi \right) \ln \left( 1-\mu -\xi \right) -\left( 1-\mu
^{\prime }-\xi \right) \ln \left( 1-\mu ^{\prime }-\xi \right) +  \notag \\
&&\left( \mu +\mu ^{\prime }+2\xi \right) \ln \cos \frac{1}{2}\beta +\left(
2-\mu -\mu ^{\prime }-2\xi \right) \ln \sin \frac{1}{2}\beta .  \notag
\end{eqnarray}

Relation between $\xi $ and other parameters for $n\rightarrow \infty $
follows from Eqs (\ref{b9})-(\ref{b11}), 
\begin{eqnarray}
\xi &=&\frac{2-\left( 1-t^{2}\right) \left( \mu +\mu ^{\prime }\right) -%
\sqrt{\Delta }}{2\left( 1-t^{2}\right) }+O\left( n^{-1}\right) ,\quad \text{ 
}t\neq 1,  \label{ab2} \\
\xi &=&\frac{1}{2}\left( 1-\mu ^{\prime }\right) \left( 1-\mu \right)
+O\left( n^{-1}\right) ,\quad t=1,  \notag
\end{eqnarray}
where the discriminant can be written in the form 
\begin{equation}
\Delta =\left( 2-\left( 1-t^{2}\right) \left( \mu +\mu ^{\prime }\right)
\right) ^{2}-4\left( 1-t^{2}\right) \left( 1-\mu \right) \left( 1-\mu
^{\prime }\right) .  \label{ab2.1}
\end{equation}
Since $\Delta \geqslant 0$ and $4\left( 1-t^{2}\right) \left( 1-\mu \right)
\left( 1-\mu ^{\prime }\right) \geqslant 0$ this results that the principal
term $\xi \geqslant 0.$ The residual $O\left( n^{-1}\right) $ does not
affect bound as the asymptotic constant can be corrected, while the
principal term can be used in $\lambda $ in Eq. (\ref{ab10}). So, $\lambda $
then is a function of three parameters, $\lambda =\lambda \left( \mu ,\mu
^{\prime };\beta \right) .$

\begin{figure}[htb]
\center  \includegraphics[width=0.9\textwidth]{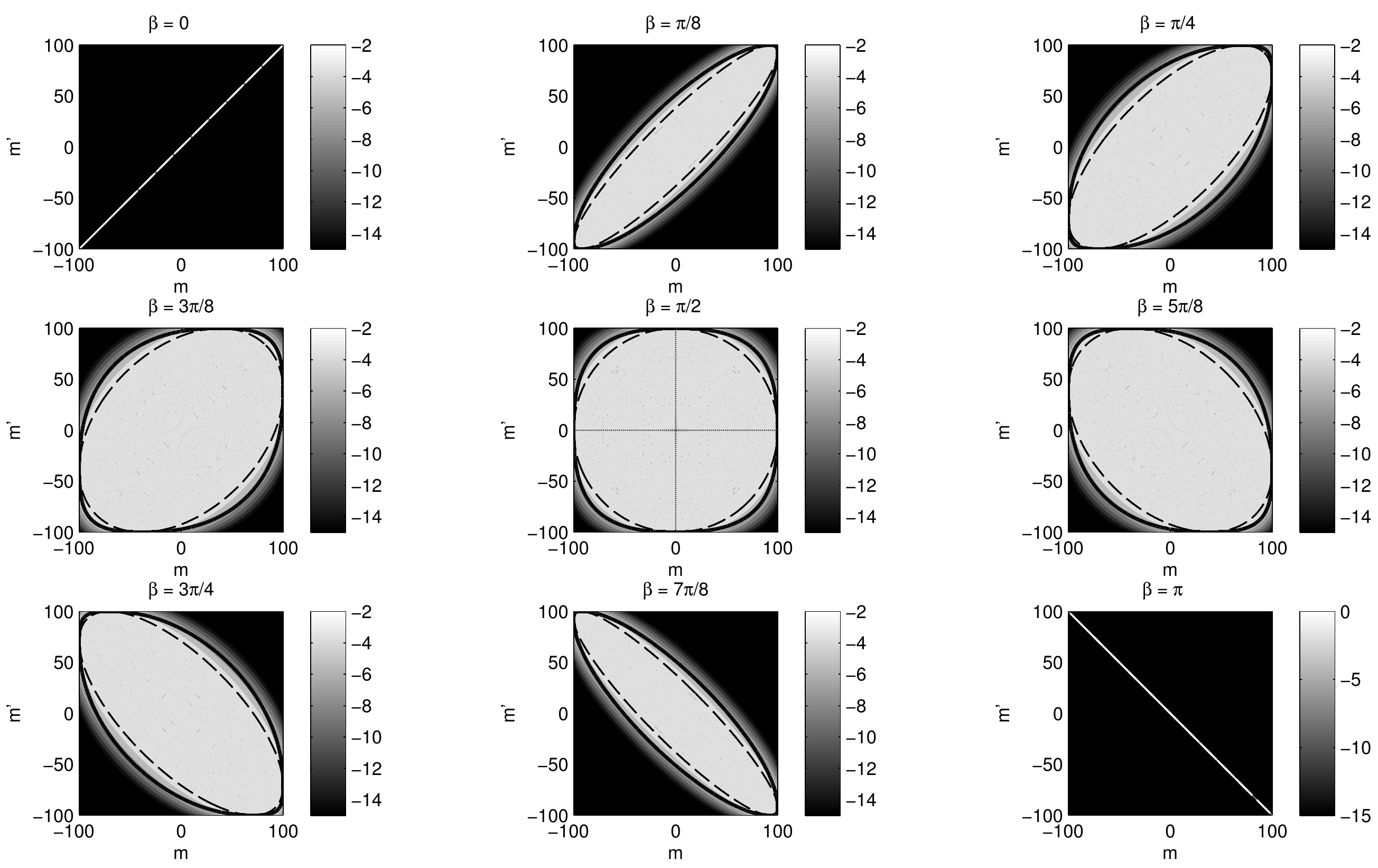}
\caption{Magnitude of rotation coefficients, $\log _{10}\left|
H_{n}^{m^{\prime }m}\left( \protect\beta \right) \right| $ at $n=100$ and
different $\protect\beta $. The solid bold curves show analytical bounds of
exponentially decaying region determined by Eq. (\protect\ref{ab11}). The
dashed curves plot the ellipse, Eq. (\protect\ref{arf12}), obtained from
asymptotic analysis of recursions for large $n$.}
\label{Fig2}
\end{figure}

Bound (\ref{ab9}) is tighter than (\ref{b2}) when $\lambda <0$, as it shows
that for $n\rightarrow \infty $ the rotation coefficients in parameter region 
$\lambda \left( \mu ,\mu ^{\prime },\beta \right) <0$ become exponentially
small. This region of exponentially small $\left| H_{n}^{m^{\prime }m}\left(
\beta \right) \right| $ at fixed $\beta $ is bounded by curve 
\begin{equation}
\lambda \left( \mu ,\mu ^{\prime };\beta \right) =0.  \label{ab11}
\end{equation}
In Fig. \ref{Fig2} computations of $\log \left| H_{n}^{m^{\prime }m}\left(
\beta \right) \right| $ at different $\beta $ and large enough $n$ $\left(
n=100\right) $ are shown. Here also the boundary curve (\ref{ab11}) is
plotted (the curve is extended by symmetry for all $\beta $ and $\mu $ and $%
\mu ^{\prime },$ so it becomes a closed curve). It is seen, that it agrees
with the computations and indeed, $\left| H_{n}^{m^{\prime }m}\left( \beta
\right) \right| $ decays exponentially.

\section{Asymptotic behavior of recursion}

\begin{figure}[htb]
\center  \includegraphics[width=0.9\textwidth]{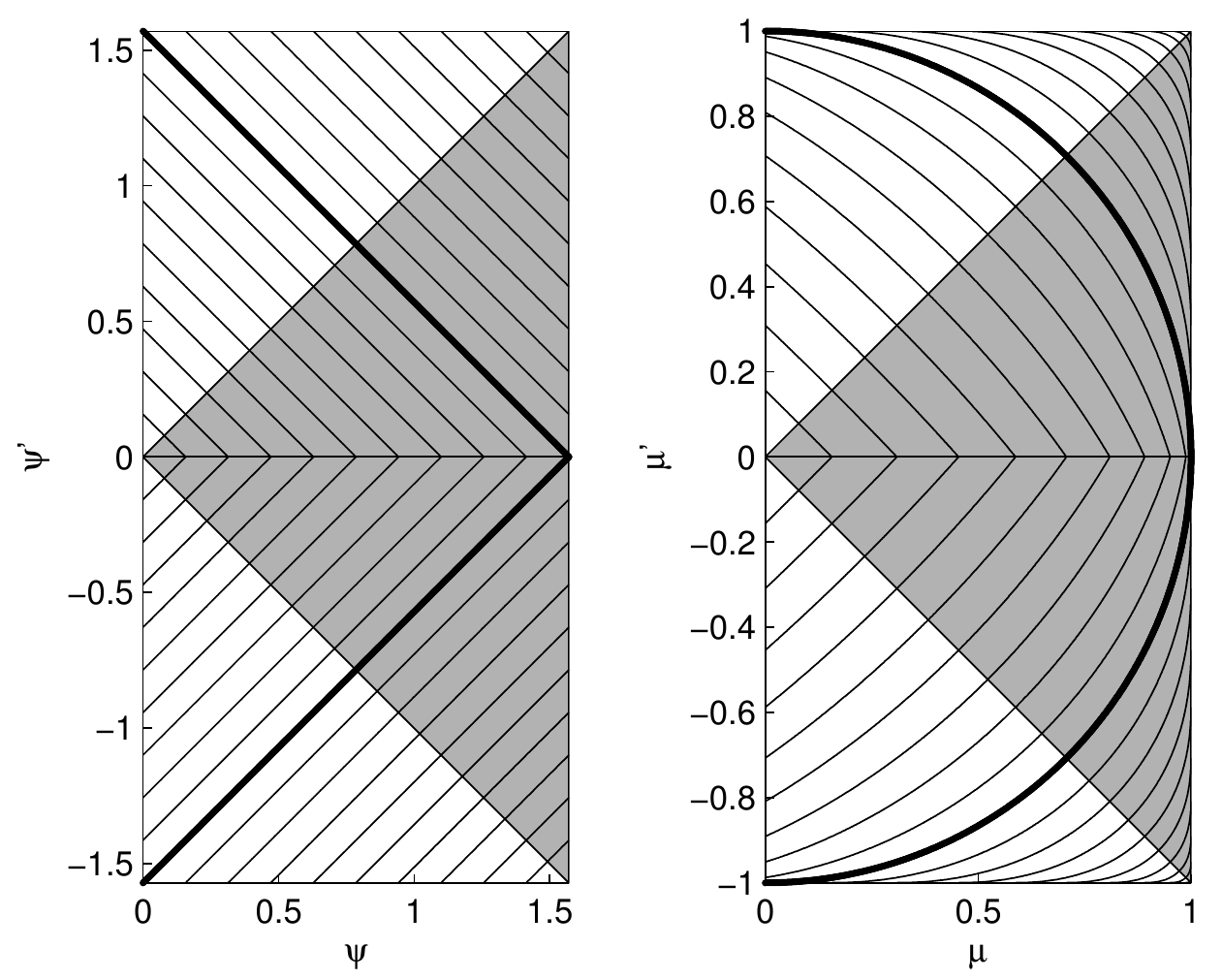}
\caption{The characteristics of equations ( \protect\ref{ar5}) and ( \protect
\ref{arf5}). The shaded area shows the region for which the rotation
coefficients should be computed as symmetry relations (\protect\ref{s1}) can
be applied to obtain the values in the other regions.}
\label{Fig3}
\end{figure}

Consider now asymptotic behavior of recursion (\ref{rec1.10}), where
coefficients of recursion $d_{n}^{m}$ do not depend on $\beta $. First, we
note that this relation can be written in the form 
\begin{gather}
k_{n}^{m^{\prime }}\left( H_{n}^{m^{\prime }+1,m}-H_{n}^{m^{\prime
}-1,m}\right) -l_{n}^{m^{\prime }}\left( H_{n}^{m^{\prime
}+1,m}+H_{n}^{m^{\prime }-1,m}\right) =  \label{arf1} \\
k_{n}^{m}\left( H_{n}^{m^{\prime },m+1}-H_{n}^{m^{\prime },m-1}\right)
-l_{n}^{m}\left( H_{n}^{m^{\prime },m-1}+H_{n}^{m^{\prime },m+1}\right) , 
\notag
\end{gather}%
where 
\begin{equation}
k_{n}^{m}=\frac{1}{2}\left( d_{n}^{m-1}+d_{n}^{m}\right) ,\quad l_{n}^{m}=%
\frac{1}{2}\left( d_{n}^{m-1}-d_{n}^{m}\right) .  \label{arf2}
\end{equation}%
At $m\neq 0$ and large $n$ and $n-\left\vert m\right\vert ,$ $m=n\mu $,
asymptotics of coefficients $d_{n}^{m}$ can be obtained from Eq. (\ref%
{rec1.11}), 
\begin{eqnarray}
k_{n}^{m} &=&\text{sgn}\left( \mu \right) \left( 1-\mu ^{2}\right) ^{1/2}%
\frac{1}{2n^{-1}}\left[ 1+O(n^{-1})\right] ,  \label{arf3} \\
l_{n}^{m} &=&\frac{1}{4}\text{sgn}\left( \mu \right) \frac{\mu }{\left(
1-\mu ^{2}\right) ^{1/2}}\left[ 1+O(n^{-1})\right] .  \notag
\end{eqnarray}%
Hence, asymptotically relation (\ref{arf1}) turns into 
\begin{eqnarray}
\!\!\!\!\!\!\!\!\!\!\!\!\!\!\!\!\!\!&&\text{sgn}\left( \mu ^{\prime }\right) \left[ 
\frac{(1-\mu ^{\prime 2})^{1/2}\left(H_{n}^{m^{\prime }+1,m}-H_{n}^{m^{\prime }-1,m}\right)}{2h}-\frac{}{}\frac{\mu
^{\prime }\left( H_{n}^{m^{\prime
}+1,m}+H_{n}^{m^{\prime }-1,m}\right)}{4\left( 1-\mu ^{\prime 2}\right) ^{1/2}}\right]   \label{arf4} \\
\!\!\!\!\!\!\!\!\!\!\!\!\!\!\!\!\!\!&=&\text{sgn}\left( \mu \right) \left[ \frac{(1-\mu ^{2})^{1/2}\left(H_{n}^{m^{\prime },m+1}-H_{n}^{m^{\prime },m-1}\right)}{2h}- \frac{\mu\left(H_{n}^{m^{\prime },m+1}+H_{n}^{m^{\prime
},m-1}\right)}{4\left(1-\mu ^{2}\right) ^{1/2}}\right],  \notag
\end{eqnarray}%
where $ h={1}/{n}$. Let us intepret now $H_{n}^{m^{\prime }m}$ as samples of differentiable
function\ $H_{n}\left( \mu ^{\prime },\mu \right) $ on a $\left( 2n+1\right)
\times \left( 2n+1\right) $ grid on the square $\left( \mu ^{\prime },\mu
\right) \in $ $\left[ -1,1\right] \times \left[ -1,1\right] $ with step $h$
in each direction, $H_{n}\left( m^{\prime }/n,m/n\right) =H_{n}^{m^{\prime
}m}$. In this case relation (\ref{arf4}) corresponds to a central
difference scheme for the hyperbolic PDE 
\begin{eqnarray}
&&\text{sgn}\left( \mu \right) \left( 1-\mu ^{2}\right) ^{1/2}\frac{\partial
H_{n}}{\partial \mu }-\text{sgn}\left( \mu ^{\prime }\right) \left( 1-\mu
^{\prime 2}\right) ^{1/2}\frac{\partial H_{n}}{\partial \mu ^{\prime }}-  \label{arf5}\\
&&\frac{1}{2}\left[ \text{sgn}\left( \mu \right) \frac{\mu }{\left( 1-\mu
^{2}\right) ^{1/2}}-\text{sgn}\left( \mu ^{\prime }\right) \frac{\mu
^{\prime }}{\left( 1-\mu ^{\prime 2}\right) ^{1/2}}\right] H_{n}=0, \notag
\end{eqnarray}%
where  $H_{n}$ is approximated to $O(h)$ via its values at neighbouring grid points in
each direction, 
\begin{eqnarray}
&&H_{n}\left( \mu ^{\prime },\mu \right) =\frac{1}{2}\left(
H_{n}\left( \mu ^{\prime }-h,\mu \right) +H_{n}\left( \mu ^{\prime }+h,\mu
\right) +O\left( h\right) \right) = \notag \\
&& \frac{1}{2}\left( H_{n}\left( \mu
^{\prime },\mu -h\right) +H_{n}\left( \mu ^{\prime },\mu +h\right) +O\left(
h\right) \right).
\end{eqnarray}
Note then  $K_{n}\left( \mu ^{\prime },\mu \right) $ defined as 
\begin{equation}
K_{n}\left( \mu ^{\prime },\mu \right) =\left( 1-\mu ^{\prime 2}\right)
^{1/4}\left( 1-\mu ^{2}\right) ^{1/4}H_{n}\left( \mu ^{\prime },\mu \right) ,
\label{arf6}
\end{equation}
satisfies 
\begin{equation}
\text{sgn}\left( \mu \right) \left( 1-\mu ^{2}\right) ^{1/2}\frac{\partial
K_{n}}{\partial \mu }-\text{sgn}\left( \mu ^{\prime }\right) \left( 1-\mu
^{\prime 2}\right) ^{1/2}\frac{\partial K_{n}}{\partial \mu ^{\prime }}=0.
\label{arf7}
\end{equation}
Let us introduce the variables $\psi =\arcsin \mu$ and  $\psi ^{\prime }=\arcsin \mu ^{\prime}$, and 
\begin{equation}
G_{n}( \psi ^{\prime },\psi) =K_{n}( \mu ^{\prime },\mu) ,\quad -\frac{\pi}{2}\!\leqslant\! \psi ,\psi ^{\prime }\!\leqslant\! \frac{\pi}{2}.\label{ar4}
\end{equation}
In this case $\mu =\sin \psi ,$ $(1-\mu ^{2})^{1/2}=\cos \psi $ and
similarly, $\mu ^{\prime }=\sin \psi ^{\prime },$ $(1-\mu ^{\prime
2})^{1/2}=\cos \psi ^{\prime }$. So in these variables Eq. (\ref{arf3}) turns
into 
\begin{equation}
\text{sgn}\left( \psi ^{\prime }\right) \frac{\partial G_{n}}{\partial \psi
^{\prime }}-\text{sgn}\left( \psi \right) \frac{\partial G_{n}}{\partial
\psi }=0.  \label{ar5}
\end{equation}
Due to symmetry relations (\ref{s1}), we can always constrain the region
with $\mu \geqslant 0$ $\left( 0\leqslant \psi \leqslant \pi /2\right) $, so
sgn$\left( \psi \right) =1$. So, we have two families of characteristics of
this equation in this region, 
\begin{equation}
\psi ^{\prime }+\psi =C^{+},\quad 0\leqslant \psi ^{\prime }\leqslant \frac{\pi}{2};  \quad
\psi ^{\prime }-\psi =C^{-},\quad -\frac{\pi}{2}\leqslant \psi ^{\prime
}\leqslant 0,  \label{arf8}
\end{equation}
which are also characteristics of equation (\ref{arf5}) for $\psi =\arcsin
\mu ,$ $\psi ^{\prime }=\arcsin \mu ^{\prime }$. Figure \ref{Fig3}
illustrates these characteristics in the $\left( \psi ,\psi ^{\prime }\right) $
and the $\left( \mu ,\mu ^{\prime }\right) $ planes. It is interesting to
compare these characteristics with the results shown in Fig. \ref{Fig2}. One can
expect that the curves separating the regions of exponentially small values
of the rotation coefficients should follow, at least qualitatively,
 some of the characteristic curves. Indeed, while $H_{n}\left( \mu ^{\prime },\mu
\right) $ according to Eq. (\ref{arf6}) can change along the characteristics
(the value of $K_{n}\left( \mu ^{\prime },\mu \right) $ should be constant),
such a change is rather weak (proportional to $\left( 1-\mu ^{\prime 2}\right)
^{1/4}\left( 1-\mu ^{2}\right) ^{1/4}$, which cannot explain the exponential
decay). However, we can see that the boundaries of the regions plotted at
different $\beta $ partly coincide with the characteristics (e.g. at $\beta
=\pi /4$ this curve qualitatively close to characteristic family $C^{-}$ at $%
\mu ^{\prime }<0$, but qualitatively different from the characteristic
family $C^{+}$ for $\mu ^{\prime }>0$; similarly the curve $\beta =3\pi /4$
 coincides with one of characteristics $C^{+}$ for $\mu ^{\prime }>0$,
while characteristics of family $C^{-}$ are rather orthogonal to the curve
at $\mu ^{\prime }<0$). As Eq. (\ref{ar5}) should be valid for any $\beta $
this creates a puzzle, the solution of which can be explained as follows. 

\begin{figure}[htb]
\center  \includegraphics[width=0.9\textwidth]{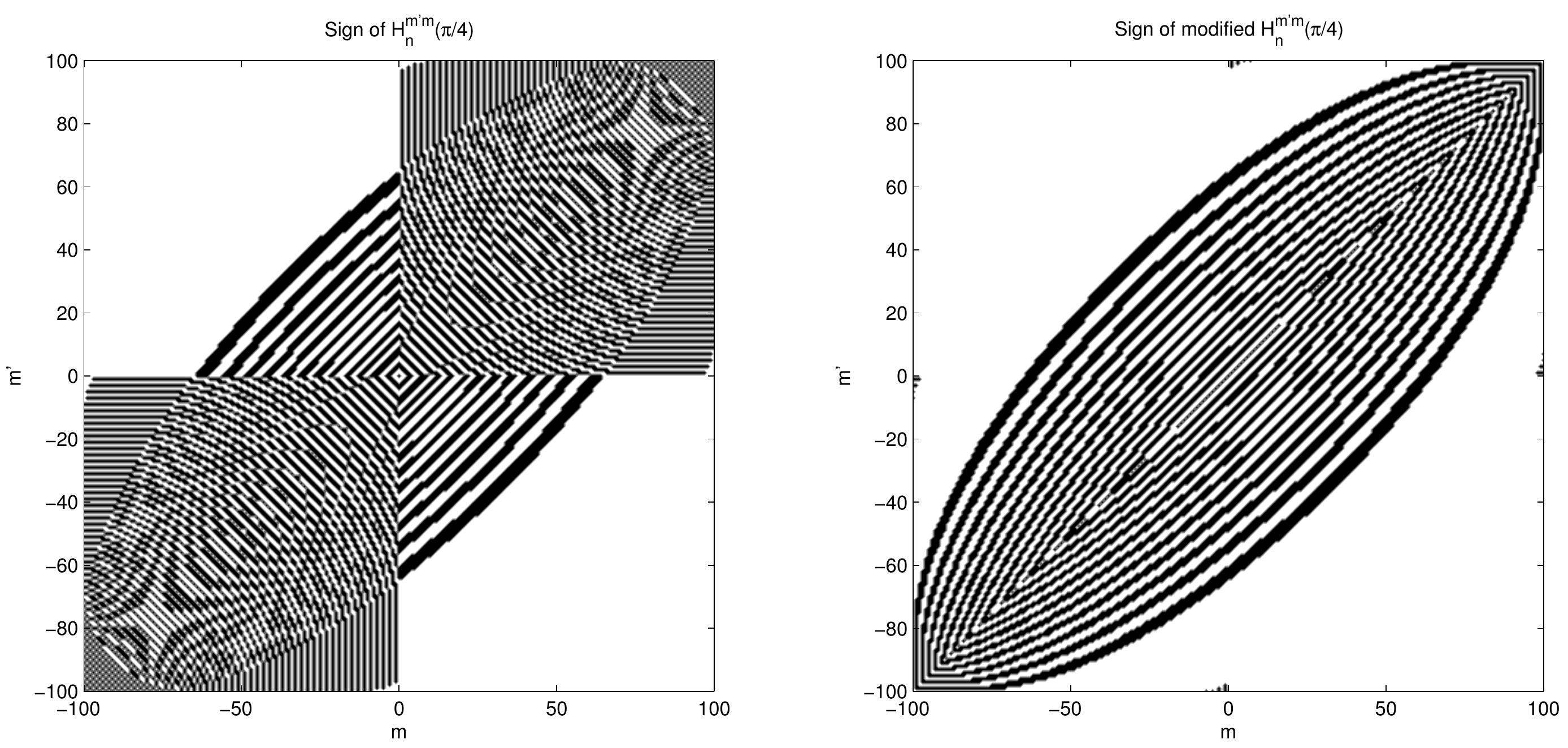}
\caption{The sign of coefficients $H_{n}^{m^{\prime }m}$ and modified
coefficients $\widehat{H}_{n}^{m^{\prime }m}$, Eq. (\protect\ref{arf9}), at $%
\protect\beta =\protect\pi /4$. Positive sign (white) is assigned to coefficients which magnitude is below $10^{-13}$.}
\label{Fig4}
\end{figure}

 $H_{n}^{m^{\prime }m}\left( \beta \right) $ considered
as a function of $m^{\prime }$ oscillates with some local frequency $\omega
_{n}\left( \mu ^{\prime },\mu ,\beta \right) $. In the regions where this
frequency is small, the transition from the discrete relation (\ref{arf4}) to the
PDE (\ref{arf5}) is justified. However at frequencies $\omega _{n}\left(
\mu ^{\prime },\mu ,\beta \right) \sim 1$ the PDE is not valid. Such regions exist, e.g., if $H_{n}^{m^{\prime }m}\left( \beta \right)
=\left( -1\right) ^{m^{\prime }}\widehat{H}_{n}^{m^{\prime }m}\left( \beta
\right) ,$ where $\widehat{H}_{n}^{m^{\prime }m}$ oscillates with a low
frequency. For example, Eqs (\ref{v3}) and (\ref{r6}) show that for $0<\beta
<\pi $ the boundary value $H_{n}^{m^{\prime }n}\left( \beta \right) $ is a 
smooth function of $\mu ^{\prime }$ for $\mu ^{\prime }<0$, while it cannot
be considered as differentiable function of $\mu ^{\prime }$ for $\mu
^{\prime }>0$. On the other hand, this example shows that the function $\widehat{%
H}_{n}^{m^{\prime }n}\left( \beta \right) $ has smooth behavior for $\mu
^{\prime }>0$ and a differential equation can be considered for this
function. Equations (\ref{arf1}) and (\ref{arf4}) show that for
function $\widehat{H}_{n}^{m^{\prime }m}$ one obtains the same equation, but the sign of sgn$\left( \mu ^{\prime }\right) $ should be changed.
Particularly, this means that if this is the case then characteristics of the
family $C^{+}$ can be extended to the region $-\pi /2\leqslant \psi ^{\prime }\leqslant 0,$ while the characteristics of the family $C^{-}$ can be continued to
the region $0\leqslant \psi ^{\prime }\leqslant \pi /2$. Of course, such an
extension must be done carefully, based on the analysis, and 
 this also depends on the values of $\beta $, which plays the role of a parameter. Fig. \ref{Fig4} illustrates the signs of function $H_{n}^{m^{\prime }m}\left( \beta \right) $ and 
\begin{equation}
\widehat{H}_{n}^{m^{\prime }m}\left( \beta \right) =\left\{ 
\begin{array}{c}
\epsilon _{m^{\prime }}\epsilon _{-m}H_{n}^{m^{\prime }m}\left( \beta
\right) ,\quad m<m^{\prime } \\ 
\epsilon _{-m^{\prime }}\epsilon _{m}H_{n}^{m^{\prime }m}\left( \beta
\right) ,\quad m\geqslant m^{\prime }%
\end{array}
\right. .  \label{arf9}
\end{equation}
It is seen that $\widehat{H}_{n}^{m^{\prime }m}$ is a ``smoother'' function
of $m^{\prime }$ and $m$ (for $0<\beta <\pi /2$; this is not the case for $%
\pi /2<\beta <\pi $; for those values we use the third symmetry relation (\ref{s1}). Also note that for $m\geqslant m^{\prime }$ the function $\widehat{H}_{n}^{m^{\prime }m}$ coincides with $d_{n}^{m^{\prime}m}$ due to relation (\ref{r8}).

This enables determination of the boundary curve separating
oscillatory and exponentially decaying regions of $\widehat{H}%
_{n}^{m^{\prime }m}$. Indeed, at $m^{\prime }=0$, $m\geqslant 0$ the
boundary value (\ref{v2}) and the first symmetry (\ref{s1}) provides that $%
\widehat{H}_{n}^{0m}\left( \beta \right) $ is proportional to the
associated Legendre function $P_{n}^{m}(\cos \beta )$. Function $%
P_{n}^{m}(x) $ satisfies differential equation 
\begin{equation}
\left( 1-x^{2}\right) \frac{d^{2}w}{dx^{2}}-2x\frac{dw}{dx}+\left[ n(n+1)-%
\frac{m^{2}}{1-x^{2}}\right] w=0,  \label{arf10}
\end{equation}
which for  
\begin{equation}
y(x)=\left( 1-x^{2}\right) ^{1/2}w(x),  \label{arf10.1}
\end{equation}
at large $n$ and $m=\mu n$ turns into 
\begin{equation}
\frac{d^{2}y}{dx^{2}}=n^{2}q\left( x\right) y,\quad q\left( x\right) =-\frac{%
1-x^{2}-\mu ^{2}}{\left( 1-x^{2}\right) ^{2}},  \label{arf10.2}
\end{equation}

An accurate asymptotic study can be done based on the Liouville-Green or WKB-approximation \cite{Naifeh73:Book}, while here we limit ourselves with the
qualitative observation, that $q(x_{\mu })=0$ at $1-x_{\mu }^{2}-\mu ^{2}=0$%
, which is a ``turning'' point, such that at $\mu ^{2}>1-x_{\mu }^{2}$ we
have $q\left( x\right) >0$ which corresponds to asymptotically
growing/decaying regions of $y$ (the decaying solution corresponds to the associated Legendre functions of the first kind, which is our case). Region $%
\mu ^{2}<1-x_{\mu }^{2}$ corresponds to $q\left( x\right) <0$ and to the
oscillatory region. The vicinity of the turning point can be studied separately (using the Airy functions) \cite{Naifeh73:Book}, but it should be
noticed immediately that the local frequency $\omega \sim n\sqrt{-q}$ is
much smaller than $n$ at $\left| q\right| \ll 1$, so function $\widehat{H}%
_{n}^{m^{\prime }m}$ is relatively smooth on the grid with step $h=n^{-1}$ in the vicinity of this turning point. Hence, for the characteristic $C^{-}$
passing through the turning point $\mu =\sqrt{1-x_{\mu }^{2}}=\sqrt{1-\cos
^{2}\beta }=\sin \beta $ at $\mu ^{\prime }=0$ ($\psi ^{\prime }=0$) 
\begin{equation}
C^{-}=\left. \psi ^{\prime }-\psi \right| _{\psi ^{\prime }=0}=-\psi
=-\arcsin \mu =-\beta ,\quad 0\leqslant \beta \leqslant \pi /2.
\label{arf11}
\end{equation}
(similarly, characteristic $C^{+}$ can be considered, which provides the
boundary curve for the case $\pi /2\leqslant \beta \leqslant \pi ).$ We note
now that curve $\psi ^{\prime }-\psi =-\beta $ in $\left( \mu ,\mu ^{\prime
}\right) $ space describes a piece of ellipse. Using symmetries (\ref{s1})
we can write equation for this ellipse in the form 
\begin{equation}
\frac{\left( \mu +\mu ^{\prime }\right) ^{2}}{4\cos ^{2}\frac{1}{2}\beta }+%
\frac{\left( \mu -\mu ^{\prime }\right) ^{2}}{4\sin ^{2}\frac{1}{2}\beta }=1.
\label{arf12}
\end{equation}
The ellipse has semiaxes $\cos \frac{1}{2}\beta $ and $\sin \frac{1}{2}\beta 
$, which are turned to $\pi /4$ angles in the $\left( \mu ,\mu ^{\prime
}\right) $. This ellipse is also shown in Fig. \ref{Fig2}, and it is seen that it approximates the regions of decay obtained from the analysis of bounds of the rotation coefficients. Also the ellipse is located completely inside those bounds and provides somehow tighter (but not strictly proven) boundaries.

More accurate consideration and asymptotic behavior of the rotation coefficients at large $n$ can be obtained using the PDE and the boundary values of the coefficients (\ref{v2}) and (\ref{v3}). Such analysis, however, deserves a separate paper and is not presented here, as the present goal is to provide a qualitative picture and develop a stable numerical procedure.

\section{Stability of recursions}
Now we consider stability of recursion (\ref{rec1.10}), which can be used to determine the rotation coefficients at $m\geqslant 0$ and $\left| m^{\prime
}\right| \leqslant m$ for $0\leqslant $ $\beta \leqslant \pi $ as for all
other values of $m$ and $m^{\prime }$ symmetries (\ref{s1}) can be used. This recursion is two dimensional and, in principle, one can resolve it with respect to any of its terms, e.g. $H_{n}^{m^{\prime }+1,m}$, and propagate it in the direction of increasing $m^{\prime }$ if the initial and boundary values are known. Several steps of the recursion can be performed anyway and there is no stability question for relatively small $n$. However, at large $n $ stability becomes critical, and, so the asymptotic analysis and behavior plays an important role for establishing of stability conditions.

\subsection{Courant-Friedrichs-Lewy (CFL) condition}
Without any regard to a finite difference approximation of a PDE recursion (\ref{rec1.10}) can be written in the form (\ref{arf1}), which for large $n$ takes
the form (\ref{arf4}). The principal term of (\ref{arf1}) for $n\rightarrow \infty $ here can be written as 
\begin{equation}
H_{n}^{m^{\prime }+1,m}-H_{n}^{m^{\prime }-1,m}=c\left( H_{n}^{m^{\prime
},m+1}-H_{n}^{m^{\prime },m-1}\right) ,\quad c=\frac{k_{n}^{m}}{%
k_{n}^{m^{\prime }}}.  \label{CFL1}
\end{equation}
An analysis of a similar recursion, appearing from the two-wave equation is provided in \cite{Courant28:MathAn}, which can be also applied to the one-wave equation approximated by the central difference scheme. If we treat here $m^{\prime }$ as an analog of time, $m$ as an analog of a spatial variable, and $c$ as the wave speed (the grid in both variables has the same
step $\Delta m=\Delta m^{\prime }=1$), then the Courant-Friedrichs-Lewy (CFL) stability condition becomes 
\begin{equation}
\left| c\right| \leqslant 1.  \label{CFL2}
\end{equation}
Note now that from definitions (\ref{arf2}) and (\ref{rec1.11}) we have 
\begin{equation}
k_{n}^{0}=0,\quad k_{n}^{-m}=-k_{n}^{m},\quad k_{n}^{m}\geqslant
k_{n}^{m+1},\quad m=1,...,n-1.  \label{CFL3}
\end{equation}
The CFL condition is satisfied for any $m\geqslant \left| m^{\prime
}\right| $, $m^{\prime }\neq 0.$ This is also consistent with the asymptotic
behavior of the recursion coefficients (\ref{arf3}), as 
\begin{equation}
c^{2}\sim \frac{1-\mu ^{2}}{1-\mu ^{\prime 2}},  \label{CFL4}
\end{equation}
which shows that in region $\mu ^{\prime 2}\leqslant \mu ^{2}$ we have Eq. (\ref{CFL2}). Note that the CFL condition for the central difference
scheme includes only the absolute value of $c$ so, independent of  which variable $H_{n}^{m^{\prime }+1,m}$ or $H_{n}^{m^{\prime
}-1,m}$ the recursion (\ref{CFL1}) is resolved the scheme satisfies the
necessary stability condition (CFL). This means that
within the region $m\geqslant \left| m^{\prime }\right| $ the scheme can be
applied in the forward or backward directions, while some care may be needed for passing the value $m^{\prime }=0$. 
\begin{figure}[htb]
\center  \includegraphics[width=0.9\textwidth]{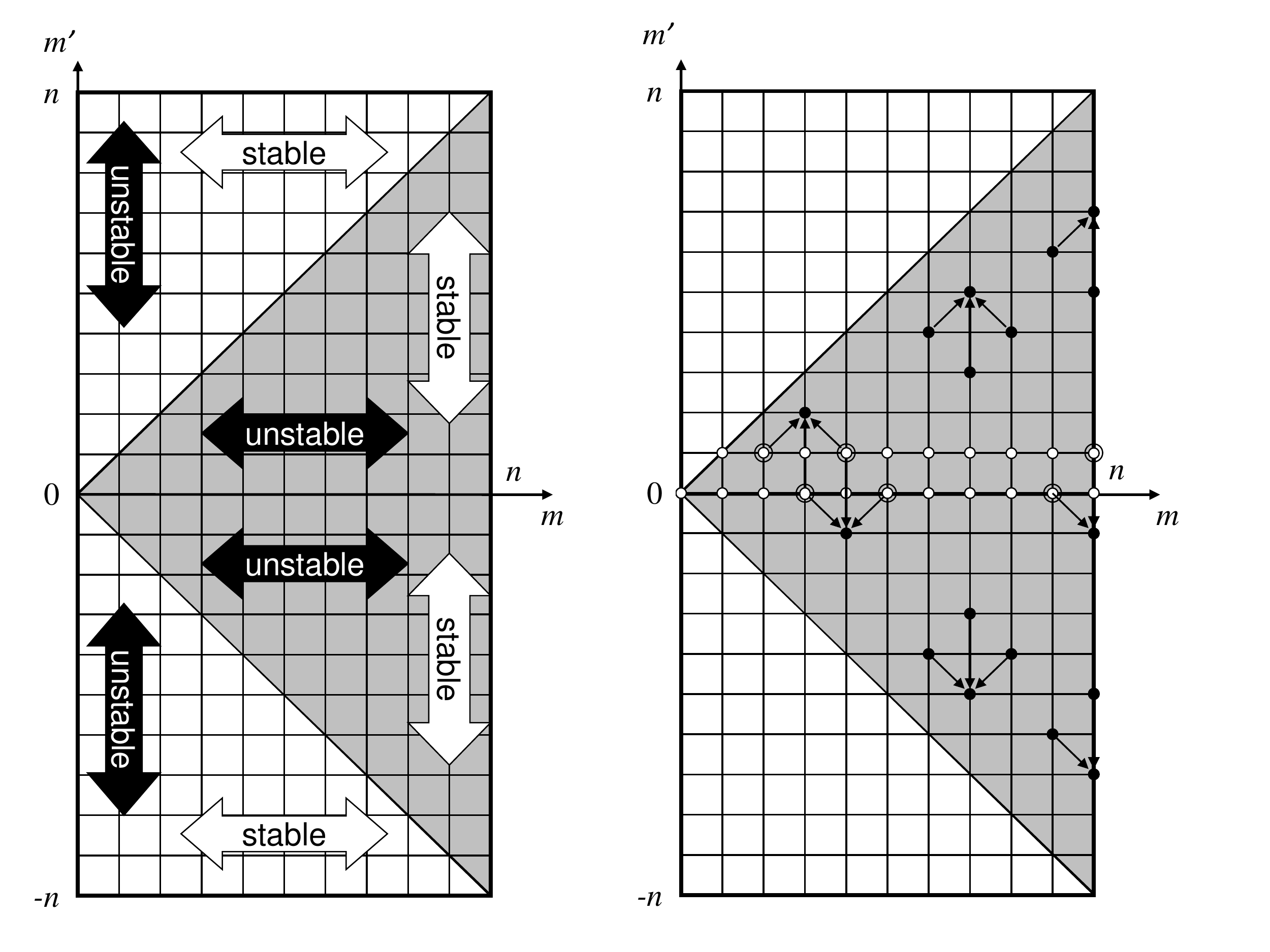}
\caption{On the left are shown unstable and conditionally stable
directions of propagation for the recursion (\protect\ref{CFL1}) based on the
Courant-Friedrichs-Lewy (CFL) criterion. On the right are shown 
the stencils for the recursive algorithm  in the shaded
region. The nodes with white discs are the initial values of the recursion
computed using Eqs (\protect\ref{v2}) and (\protect\ref{a1}).}
\label{Fig5}
\end{figure}

This analysis shows also that if recursion (\ref{CFL1}) will be resolved
with respect to $H_{n}^{m^{\prime },m+1}$ or $H_{n}^{m^{\prime },m-1}$ in
region $m\geqslant \left| m^{\prime }\right| $ then the recursion will be
absolutely unstable, as it does not satisfy the necessary condition, which
in this case will be $\left| 1/c\right| \leqslant 1$ (as the recursion is
symmetric with respect to $m$ and $m^{\prime }$). Figure \ref{Fig5} (chart
on the left) shows the stable and unstable directions of propagation. Under
``stable'' we mean here conditional, or neutral, stability, as the CFL
criterium is only a necessary, not sufficient, condition.

\subsection{Von Neumann stability analysis}
The von Neumann, or Fourier, stability analysis is a usual tool for investigation of finite difference schemes of linear equations with constant coefficients (original publication \cite{Charney50:QJG}, various applications can be found elsewhere). Despite the recursion we study is linear it  has variable coefficients. So, we can speculate that
only in some region, where such variability can be neglected, and we can
perform several recursive steps with quasiconstant coefficients, can such an analysis give us some insight on the overall stability. As the recursion
(\ref{rec1.10}) can be written in the form (\ref{arf1}), the asymptotic behavior
of the recursion coefficients (\ref{arf3}) shows that the assumption that
these coefficients are quasiconstant indeed is possible in a sense that many grid points can be handled with the same value of coefficients, as they are functions of ``slow'' variables $\mu $ and $\mu ^{\prime }$, but regions $\mu \rightarrow 1$ and $\mu ^{\prime }\rightarrow 1$ are not treatable with
this approach, as either the coefficients or their derivatives become unbounded. Hence, we apply the fon Neumann analysis for $1-\mu $ and $1-\mu
^{\prime }$ treated as quantities of the order of the unity.

Equation (\ref{arf4}) then can be written in the form 
\begin{eqnarray}
&& H_{n}^{m^{\prime }+1,m}-H_{n}^{m^{\prime }-1,m}-\frac{a}{n}\left(
H_{n}^{m^{\prime }+1,m}+H_{n}^{m^{\prime }-1,m}\right) = \label{Neu1} \\ && c\left(
H_{n}^{m^{\prime },m+1}-H_{n}^{m^{\prime },m-1}\right) -\frac{b}{n}\left(
H_{n}^{m^{\prime },m+1}+H_{n}^{m^{\prime },m-1}\right) ,  \notag
\end{eqnarray}
where $a,b,$ and $c$ are coefficients depending on $\mu $ and $\mu ^{\prime
} $ and, formally, not depending on $m$ and $m^{\prime }$ (separation to
``slow'' and ``fast'' variables typical for multiscale analysis), for $ \mu \geqslant 0$
\begin{equation}
\!\!\!\!\!\!\!\!\!a=\frac{\mu ^{\prime }}{2\left( 1-\mu ^{\prime 2}\right) },\!\quad\! b=\frac{\text{sgn}\left( \mu ^{\prime }\right) \mu }{2\left( 1-\mu ^{\prime
2}\right) ^{1/2}\left( 1-\mu ^{2}\right) ^{1/2}},\!\quad\! c=\text{sgn}\left(\mu ^{\prime }\right) \left( \frac{1-\mu ^{2}}{1-\mu ^{\prime 2}}\right)
^{1/2}.  \label{Neu1.1}
\end{equation}
Now we consider perturbation $\eta _{n}^{m^{\prime },m}$ of coefficients $%
H_{n}^{m^{\prime },m}$ and their propagation within the conditionally stable
scheme above. As the true values of\ $H_{n}^{m^{\prime },m}$
satisfy Eq. (\ref{Neu1}), the perturbation satisfies the same equation. Let $%
\widehat{\eta }_{n}^{m^{\prime }}\left( k\right) $ be the Fourier transform of $%
\eta _{n}^{m^{\prime },m}$ with respect to $m$ at layer $m^{\prime }$, where 
$k$ is the wavenumber (the $k$th harmonic of $\eta _{n}^{m^{\prime },m}$ is $%
\widehat{\eta }_{n}^{m^{\prime }}\left( k\right) e^{ikm}).$ Then Eq. (\ref%
{Neu1}) for the $k$-th harmonic takes the form 
\begin{equation}
\widehat{\eta }_{n}^{m^{\prime }+1}-\widehat{\eta }_{n}^{m^{\prime }-1}-%
\frac{b}{n}\left( \widehat{\eta }_{n}^{m^{\prime }+1}+\widehat{\eta }%
_{n}^{m^{\prime }-1}\right) =\left( 2ic\sin k-2\frac{b}{n}\cos k\right) 
\widehat{\eta }_{n}^{m^{\prime }}.  \label{Neu2}
\end{equation}
This is a one-dimensional recurrence relation, with well established  stability analysis. Particularly, one can consider solutions of type $%
\widehat{\eta }_{n}^{m^{\prime }}=\left( \lambda _{n}\right) ^{m^{\prime }}$%
, which after insertion into Eq. (\ref{Neu2}) results in the characteristic
equation 
\begin{equation}
\left( 1-\frac{a}{n}\right) \lambda _{n}^{2}-\left( 2ic\sin k-2\frac{b}{n}%
\cos k\right) \lambda _{n}-\left( 1+\frac{a}{n}\right) =0,  \label{Neu3}
\end{equation}
with roots 
\begin{equation}
\!\!\!\!\!\!\!\lambda _{n}^{\pm }\!=\!\left( 1-\frac{a}{n}\right) ^{-1}\left[ \left( ic\sin k-%
\frac{b}{n}\cos k\right)\! \pm \sqrt{\left( ic\sin k-\frac{b}{n}\cos k\right)
^{2}\!+\!1\!-\!\left( \frac{a}{n}\right) ^{2}}\right] .  \label{Neu4}
\end{equation}
If $\left| c\sin k\right| <1$ then for $n\rightarrow \infty $ we have 
\begin{eqnarray}
\left| \lambda _{n}^{\pm }\right| &\sim &\left( 1+\frac{a}{n}\right) \left[
1\mp \frac{2b}{n}\cos k\left( 1-c^{2}\sin ^{2}k+\frac{c^{2}\sin ^{2}k\cos k}{%
\sqrt{1-c^{2}\sin ^{2}k}}\right) \right] ^{1/2}  \label{Neu5} \\
&\lesssim &\left( 1+\frac{1}{n}\left[ a+\left| b\right| \left( 1+\frac{c^{2}%
}{4\sqrt{1-c^{2}}}\right) \right] \right) .  \notag
\end{eqnarray}
In the asymptotic region near $\left| c\sin k\right| =1$, since  $\left| c\right| \leqslant 1$, we also have $\left| \cos k\right| \ll 1.$ So, denoting $\left| c\sin k\right| =1-n^{-1}c^{\prime }$ and expanding $\lambda_{n}^{\pm }$ from Eq. (\ref{Neu4}) at $n\rightarrow \infty $ we obtain 
\begin{equation}
\left| \lambda _{n}^{\pm }\right| \sim \left( 1+\frac{1}{n}\left(
a+c^{\prime }\right) \right) .  \label{Neu6}
\end{equation}
Note then that for certain $k$ the recursion appears to be unstable, since $%
\left| a\right| \leqslant \left| b\right| $ in region $\left| \mu ^{\prime
}\right| <\mu .$ However, in both cases described by Eq. (\ref{Neu5}) and (%
\ref{Neu6}) the growth rate is close to one. So if we have some initial
perturbation of magnitude $\epsilon _{0}$ we have after $n$ steps (which is
the maximum number of steps for propagation from $m^{\prime }=0$ to $%
m^{\prime }=n$) error $\epsilon $  satisfies
\begin{equation}
\epsilon \leqslant \epsilon _{0}\left| \lambda _{n}\right| ^{n}\sim \epsilon
_{0}\left( 1+\frac{C}{n}\right) ^{n}\sim \epsilon _{0}e^{C}.  \label{Neu7}
\end{equation}
Here $C$ is some constant of order of 1, which does not depend on $n$, so
despite the instability the error should not grow more than a certain finite
value, $\epsilon /\epsilon _{0}\lesssim e^{C}$.

\section{Algorithms for Computation of Rotation Coefficients}
We present below two different and novel algorithms for computation of the rotation coefficients $H_{n}^{m^{\prime },m}$. The recursive algorithm is more practical (faster), while the Fast Fourier Transform (FFT) based algorithm has an advantage that it does not use any recursion and so it is free from recursion related instabilities. Availability of an
alternative independent method enables cross-validation and error/performance studies.

\subsection{Recursive algorithm}
The analysis presented above allows us to propose an algorithm for computation of the rotation coefficients based on recursion (\ref{rec1.10}). Note that this recursion, in a shortened form, is also valid for the boundary points, i.e. it holds at $m=n$ where one should set $H_{n}^{m^{\prime },n+1}=0$ (this appears autmatically as also $d_{n}^{n}=0$). Using this observation one can avoid some extra work of direct computation of the boundary values (\ref{v3}), which, however, is also not critical for the overal algorithm complexity. In the algorithm coefficients $H_{n}^{m^{\prime }m}$ are computed for each subspace $n$ independently for $m$ and $m^{\prime }$ located inside a triangle, i.e. for
values $m=0,...,n,$ $m^{\prime }=-m,...,m.$ Angle $\beta $ can take any value from $0$ to $\pi $. 

\begin{algorithmic}[1]
\State If $n=0$ set $H_{0}^{00}=1$. For other $n=1,...,p-1$ consider the rest of the algorithm.

\State Compute values $H_{n}^{0,m}\left( \beta \right) $ for $m=0,...,n$ and $
H_{n+1}^{0,m}\left( \beta \right) $ for $m=0,...,n+1$ using Eq. (\ref{v2}) (one can replace there $m^{\prime }$ with $m$ due to symmetry). It is instructive to compute these values using a stable standard routine for computation of the normalized associated Legendre functions (usually based on recursions), which avoids computation of factorials of large numbers. A
standard Matlab function serves as an example of such a routine.

\State Use relation (\ref{rec1}) to compute $H_{n}^{1,m}\left( \beta \right) $, $m=1,...,n$. Using symmetry and shift of the indices this relation can be written as 
\begin{equation}
\!\!b_{n+1}^{0}H_{n}^{1,m}= 
\frac{b_{n+1}^{-m-1}( 1-\cos \beta)}{2}
H_{n+1}^{0,m+1}-\frac{b_{n+1}^{m-1}(1+\cos \beta)}{2} H_{n+1}^{0,m-1} -a_{n}^{m}\sin \beta H_{n+1}^{0,m}.  \label{a1}
\end{equation}

\State Recursively compute $H_{n}^{m^{\prime }+1,m}\left( \beta \right) $ for 
$m^{\prime }=1,...,n-1,$ $m=m^{\prime },...,n$ using relation (\ref{rec1.10}) resolved with respect to $H_{n}^{m^{\prime }+1,m}$ 
\begin{equation}
\!\!\!\!\!\!\!\!d_{n}^{m^{\prime }}H_{n}^{m^{\prime }+1,m}= d_{n}^{m^{\prime}-1}H_{n}^{m^{\prime }-1,m}-d_{n}^{m-1}H_{n}^{m^{\prime
},m-1}+d_{n}^{m}H_{n}^{m^{\prime },m+1} ,  \label{a2}
\end{equation}
which for $m=n$ turns into 
\begin{equation}
H_{n}^{m^{\prime }+1,m}=\frac{1}{d_{n}^{m^{\prime }}}\left( d_{n}^{m^{\prime
}-1}H_{n}^{m^{\prime }-1,m}-d_{n}^{m-1}H_{n}^{m^{\prime },m-1}\right).  \label{a3}
\end{equation}

\State Recursively compute $H_{n}^{m^{\prime }-1,m}\left( \beta \right) $ for 
$m^{\prime }=-1,...,-n+1,$ $m=-m^{\prime },...,n$ using relation (\ref{rec1.10}) resolved with respect to $H_{n}^{m^{\prime }-1,m}$
\begin{equation}
\!\!\!\!\!\!{d_{n}^{m^{\prime }-1}}H_{n}^{m^{\prime }-1,m}=
d_{n}^{m^{\prime }}H_{n}^{m^{\prime }+1,m}+d_{n}^{m-1}H_{n}^{m^{\prime
},m-1}-d_{n}^{m}H_{n}^{m^{\prime },m+1} ,  \label{a4}
\end{equation}
which for $m=n$ turns into 
\begin{equation}
H_{n}^{m^{\prime }-1,m}=\frac{1}{d_{n}^{m^{\prime }-1}}\left(
d_{n}^{m^{\prime }}H_{n}^{m^{\prime }+1,m}+d_{n}^{m-1}H_{n}^{m^{\prime
},m-1}\right).  \label{a5}
\end{equation}

\State Apply the first and the second symmetry relations (\ref{s1}) to obtain
all other values $H_{n}^{m^{\prime }m}$ outside the computational triangle $%
m=0,...,n,$ $m^{\prime }=-m,...,m.$
\end{algorithmic}

Figure \ref{Fig5} (right) illustrates this algorithm. 
It is clear that the algorithm needs $O(1)$ operations per value of $H_{n}^{m^{\prime },m}$. It also can be applied to each subspace independently, and is parallelizable. So, the complexity for a single subspace of degree $n$ is $O\left( n^{2}\right) $, and the cost to compute all the rotation coefficients
for $p$ subspaces ($n=0,...,p-1$) is $O\left( p^{3}\right) $. It also can be
noticed that for computation of rotation coefficients for all subspaces $n=0,...,p-1$ the algorithm can be simplified, as instead of computation of $H_{n}^{0,m}$and $H_{n+1}^{0,m}$ for each subspace in step 2, $H_{n}^{0,m}$ can be computed for all $n=1,...,p$ ($m=0,...,n$)
and stored. Then the required initial values can be retrieved at the time of processing of the $n$th subspace.

\subsection{FFT based algorithms}
\subsubsection{Basic algorithm}
We propose this algorithm based on Eq. (\ref{r4.1}), which for $\alpha =0$ and $\gamma =0$ takes the form 
\begin{eqnarray}
f_{n}^{m}\left( \widehat{\varphi };\beta ,\widehat{\theta }\right)
&=&\sum_{m^{\prime }=-n}^{n}F_{n}^{m^{\prime }m}\left( \beta ,\widehat{%
\theta }\right) e^{im^{\prime }\widehat{\varphi }},\quad  \label{fft1} \\
f_{n}^{m}\left( \widehat{\varphi };\beta ,\widehat{\theta }\right)
&=&Y_{n}^{m}\left( \theta \left( \widehat{\varphi },\beta ,\widehat{\theta }%
\right) ,\varphi \left( \widehat{\varphi },\beta ,\widehat{\theta }\right)
\right) ,  \notag \\
F_{n}^{m^{\prime }m}\left( \beta ,\widehat{\theta }\right)
&=&(-1)^{m^{\prime }}\sqrt{\frac{2n+1}{4\pi }\frac{(n-\left| m^{\prime
}\right| )!}{(n+\left| m^{\prime }\right| )!}}P_{n}^{\left| m^{\prime
}\right| }(\cos \widehat{\theta })H_{n}^{m^{\prime }m}\left( \beta \right) .
\notag
\end{eqnarray}
Here we used the definition of the spherical harmonics (\ref{3}); $\theta
\left( \widehat{\varphi },\beta ,\widehat{\theta }\right) $ and $\varphi
\left( \widehat{\varphi },\beta ,\widehat{\theta }\right) $ are determined
by the rotation transform (\ref{r2.1.3}) and (\ref{r2.2}), where we set $\alpha =0$ and $\gamma =0$. The rotation matrix $Q$ is symmetric (see Eq. (\ref{r2})) 
\begin{equation}
\left( 
\begin{array}{c}
x \\ 
y \\ 
z%
\end{array}
\right) =\left( 
\begin{array}{ccc}
-\cos \beta & 0 & \sin \beta \\ 
0 & -1 & 0 \\ 
\sin \beta & 0 & \cos \beta%
\end{array}
\right) \left( 
\begin{array}{c}
\widehat{x} \\ 
\widehat{y} \\ 
\widehat{z}%
\end{array}
\right) .  \label{fft2}
\end{equation}
Using relation between the Cartesian and spherical coordinates (\ref{1}), we
obtain 
\begin{eqnarray}
\sin \theta \cos \varphi &=&-\cos \beta \sin \widehat{\theta }\cos \widehat{%
\varphi }+\sin \beta \cos \widehat{\theta },  \label{fft3} \\
\sin \theta \sin \varphi &=&-\sin \widehat{\theta }\sin \widehat{\varphi }, 
\notag \\
\cos \theta &=&\sin \beta \sin \widehat{\theta }\cos \widehat{\varphi }+\cos
\beta \cos \widehat{\theta }.  \notag
\end{eqnarray}
This specifies functions $\varphi \left( \widehat{\varphi },\beta ,\widehat{%
\theta }\right) $ and $\theta \left( \widehat{\varphi },\beta ,\widehat{%
\theta }\right) $, $0\leqslant \varphi <2\pi ,$ $0\leqslant \theta \leqslant
\pi .$

Now, let us fix some $\widehat{\theta }$, such that $\cos \widehat{\theta }$
is not a zero of the associated Legendre function $P_{n}^{m}\left( x\right) $
at any $m=0,...,n$. Then for a given $\beta $ function $f_{n}^{m}\left( 
\widehat{\varphi };\beta ,\widehat{\theta }\right) $ is completely defined
as a function of $\widehat{\varphi }$, while $\beta ,\widehat{\theta },m,$
and $n$ play a role of parameters. The first equation shows then that this
function has a finite Fourier spectrum ($2n+1$ harmonics,). The problem then
is to find this spectrum ($F_{n}^{m^{\prime }m}$), which can be done via the
FFT, and from that determine $H_{n}^{m^{\prime }m}\left( \beta \right) $
using the last relation (\ref{fft1}). The complexity of the algorithm for
subspace $n$ is, obviously, $O\left( n^{2}\log n\right) $ and for all
subspaces $n=0,1,...,p-1$ we have complexity $O\left( p^{3}\log p\right) .$

\subsubsection{Modified algorithm}

The problem with this algorithm is that at large $n$ the associated Legendre
functions (even the normalized ones) are poorly scaled. Analysis of Eq. (%
\ref{arf10.2}) shows that to have coefficients $F_{n}^{m^{\prime }m}$ of the
order of unity parameter $\widehat{\theta }$ should be selected as close to $%
\pi /2$ as possible. On the other hand, this cannot be exactly $\pi /2$ as
in this case $P_{n}^{\left| m^{\prime }\right| }(0)=0$ for odd values of $%
n+\left| m^{\prime }\right| $. The following trick can be proposed to fix
this.

Consider two functions $g_{n}^{(1)m}\left( \widehat{\varphi };\beta \right)
=f_{n}^{m}\left( \widehat{\varphi };\beta ,\pi /2\right) $ and $%
g_{n}^{(2)m}\left( \widehat{\varphi };\beta \right) =\left. \partial
f_{n}^{m}\left( \widehat{\varphi };\beta ,\widehat{\theta }\right) /\partial 
\widehat{\theta }\right| _{\widehat{\theta }=\pi /2}$. The first function
has spectrum $\left\{ F_{n}^{m^{\prime }m}\left( \beta ,\pi /2\right)
\right\} $, while the second function $\left\{ \left. \partial
F_{n}^{m^{\prime }m}\left( \beta ,\pi /2\right) /\partial \widehat{\theta }%
\right| _{\widehat{\theta }=\pi /2}\right\} $. Note that
$P_{n}^{\left| m^{\prime }\right| }(\cos \widehat{\theta })$ is an
even function of $\widehat{x}=\cos \widehat{\theta }$ for even $n+\left|
m^{\prime }\right| $, and an odd function of $\widehat{x}=\cos 
\widehat{\theta }$ for odd values of $n+\left| m^{\prime }\right| $. In the
latter case $x=0$ is a single zero and $\left. \partial F_{n}^{m^{\prime
}m}\left( \beta ,\pi /2\right) /\partial \widehat{\theta }\right| _{\widehat{%
\theta }=\pi /2}$ is not zero for odd $n+\left| m^{\prime }\right| $ (its
absolute value reaches the maximum at $\widehat{x}=0$), while it is zero for
even $n+\left| m^{\prime }\right| $. Hence, 
\begin{equation}
\!\!\!\!\!\!\!\!\!\!g_{n}^{m}\left( \widehat{\varphi };\beta \right) =g_{n}^{(1)m}\left( 
\widehat{\varphi };\beta \right) +g_{n}^{(2)m}\left( \widehat{\varphi }%
;\beta \right)\! =\! \left[ f_{n}^{m}( \widehat{\varphi };\beta ,%
\widehat{\theta }) +b_{n}^{m}\frac{\partial }{\partial \widehat{\theta 
}}f_{n}^{m}( \widehat{\varphi };\beta ,\widehat{\theta }) \right]_{\widehat{\theta }=\pi /2},  \label{fft4}
\end{equation}
where $\gamma _{n}^{m}\neq 0$ is an arbitrary number, has Fourier spectrum 
\begin{equation}
G_{n}^{m^{\prime }m}\left( \beta \right) =H_{n}^{m^{\prime }m}\left( \beta
\right) K_{n}^{m^{\prime }m},  \label{fft5} 
\end{equation}
where, for $m^{\prime }=2k-n,\quad k=0,...,n,$ 
$$
K_{n}^{m^{\prime }m} =(-1)^{m^{\prime }}\sqrt{\frac{2n+1}{4\pi }\frac{%
(n-\left| m^{\prime }\right| )!}{(n+\left| m^{\prime }\right| )!}} P_{n}^{\left| m^{\prime }\right| }(0), 
$$
while, for $\quad m^{\prime }=2k-n-1,\quad k=1,...,n$
$$
K_{n}^{m^{\prime }m} =-(-1)^{m^{\prime }}\sqrt{\frac{2n+1}{4\pi }\frac{%
(n-\left| m^{\prime }\right| )!}{(n+\left| m^{\prime }\right| )!}} 
\gamma _{n}^{m}(n+\left| m^{\prime }\right| )(n-\left| m^{\prime }\right|
+1)P_{n}^{\left| m^{\prime }\right| -1}(0),
$$
The values for odd $n+\left| m^{\prime }\right| $ come from the well-known
recursion for the associated Legendre functions (see \cite{Abramowitz1972}), 
\begin{equation}
\frac{d}{dx}P_{n}^{m}(x)=\frac{(n+m)(n-m+1)}{\sqrt{1-x^{2}}}%
P_{n}^{m-1}\left( x\right) +\frac{mx}{1-x^{2}}P_{n}^{m}\left( x\right)
,  \label{fft5.1}
\end{equation}
(note $P_{n}^{-1}\left( x\right) =-P_{n}^{1}\left( x\right) /(n(n+1))$),
which is evaluated at $\widehat{x}=0:$%
\begin{equation}
\!\!\!\!\!\!\!\!\!\left. \frac{d}{d\widehat{\theta }}P_{n}^{\left| m^{\prime }\right| }(\cos 
\widehat{\theta })\right| _{\widehat{\theta }=\pi /2}\!\!\!=\!-\left. \frac{d}{d%
\widehat{x}}P_{n}^{\left| m^{\prime }\right| }(\widehat{x})\right| _{%
\widehat{x}=0}\!\!\!\!\!=-(n+\left| m^{\prime }\right| )(n-\left| m^{\prime }\right|
+1)P_{n}^{\left| m^{\prime }\right| -1}\!\left( 0\right) .  \label{fft6}
\end{equation}
Note also that $P_{n}^{\left| m^{\prime }\right| }(0)$ can be simply
expressed via the gamma-function (see \cite{Abramowitz1972}) and, so $%
K_{n}^{m^{\prime }m}$ can be computed without use of the associated Legendre
functions. The magnitude of arbitrary constant $\gamma _{n}^{m}$ can be
selected based on the following observation. As coefficients $%
H_{n}^{m^{\prime }m}\left( \beta \right) $ at fixed $\beta $ large $n$
asymptotically behave as functions of $m/n$ and $m^{\prime }/n$ we can try
to have odd and even coefficients $K_{n}^{m^{\prime }m}$ and $%
K_{n}^{m^{\prime }+1,m}$ to be of the same order of magnitude. We can write
this condition and the result as 
\begin{equation}
\sqrt{\frac{(n-\left| m^{\prime }\right| )!}{(n+\left| m^{\prime }\right| )!}%
}\sim \sqrt{\frac{(n-\left| m^{\prime }\right| -1)!}{(n+\left| m^{\prime
}\right| +1)!}}\gamma _{n}^{m}(n+\left| m^{\prime }\right| +1)(n-\left|
m^{\prime }\right| ),  \label{fft7}
\end{equation}
and $\gamma _{n}^{m}\sim {1}/{n}$.
Now, we can simplify expression (\ref{fft4}) for $g_{n}^{m}$. It is
sufficient to consider only positive $m$, since for negative values we can
use symmetry (\ref{s1}), while for $m=0$ we do not need Fourier transform,
as we already have Eq. (\ref{v2}). Using definitions (\ref{fft1}) and (\ref%
{3}), we obtain 
\begin{equation}
\!\!\!\!\!\!\!\!\!\frac{\partial f_{n}^{m}}{\partial\widehat{\theta }} \!=\!(-1)^{m}\sqrt{\frac{2n+1}{4\pi }\frac{(n-m)!}{(n+m)!}}e^{im\varphi}
\left[ \frac{dP_{n}^{m}(x)}{dx}\frac{\partial \cos \theta }{\partial \widehat{\theta }}+imP_{n}^{m}(x)\frac{\partial
\varphi }{\partial \widehat{\theta }}\right] _{x=\cos \theta }. \label{fft8}
\end{equation}
Differentiating  (\ref{fft3}) w.r.t. $\widehat{\theta }$
and taking values at $\widehat{\theta }=\pi /2$
\begin{equation}
\left. \frac{\partial \cos \theta }{\partial \widehat{\theta }}\right| _{%
\widehat{\theta }=\pi /2} =-\cos \beta , \quad  \left. \frac{\partial \varphi }{\partial \widehat{\theta }}\right| _{%
\widehat{\theta }=\pi /2} =\frac{\sin \beta \sin \varphi }{\sin \theta }. \label{fft9} 
\end{equation}
We also have from relations (\ref{fft3}) at $\widehat{\theta }=\pi /2$%
\begin{equation}
x=\cos \theta =\sin \beta \cos \widehat{\varphi },\quad \cos \varphi =-\frac{%
\cos \beta \cos \widehat{\varphi }}{\sqrt{1-x^{2}}},\quad \sin \varphi =-%
\frac{\sin \widehat{\varphi }}{\sqrt{1-x^{2}}}.  \label{fft10}
\end{equation}
Using these relations and identity (\ref{fft5.1}), we can write 
\begin{eqnarray}
&&\left[ \frac{dP_{n}^{m}(x)}{dx}\frac{\partial \cos \theta }{\partial 
\widehat{\theta }}+imP_{n}^{m}(x)\frac{\partial \varphi }{\partial \widehat{%
\theta }}\right] _{x=\cos \theta }= \label{fft11}\\
&&-(n+m)(n-m+1)\cos \beta \frac{%
P_{n}^{m-1}\left( x\right) }{\sqrt{1-x^{2}}}+me^{i\varphi }\sin \beta \frac{%
P_{n}^{m}(x)}{\sqrt{1-x^{2}}}.   \notag
\end{eqnarray}
Hence, function $g_{n}^{m}\left( \widehat{\varphi };\beta \right) $
introduced by Eq. (\ref{fft4}) can be written as 
\begin{eqnarray}
\!\!\!\!\!\!\!&&g_{n}^{m}\left( \widehat{\varphi };\beta \right) =(-1)^{m}\sqrt{\frac{2n+1}{%
4\pi }\frac{(n-m)!}{(n+m)!}}e^{im\varphi } \times  \label{fft12}\\
\!\!\!\!\!\!\!&& \left[ P_{n}^{m}(x)-\frac{\gamma
_{n}^{m}}{\sqrt{1-x^{2}}}\left( (n+m)(n-m+1)\cos \beta P_{n}^{m-1}\left(
x\right) -me^{i\varphi }\sin \beta P_{n}^{m}(x)\right) \right] , \notag 
\end{eqnarray}
It may appear that $x=\pm 1$ can be potentially singular, but this is not the case. Indeed, these values can be achieved only when $\beta =\pi /2$
(see the first equation (\ref{fft10})). But in this case, we can simplify
Eq. (\ref{fft12}), as we have $\cos \beta =0,$ $\sqrt{1-x^{2}}=\left| \sin 
\widehat{\varphi }\right| ,$ and so
\begin{eqnarray}
e^{i\varphi }&=&\cos \varphi +i\sin \varphi =-\frac{\cos \beta \cos \widehat{%
\varphi }}{\sqrt{1-x^{2}}}-i\frac{\sin \widehat{\varphi }}{\sqrt{1-x^{2}}}=-i%
\text{sgn}\left( \sin \widehat{\varphi }\right) , \label{fft14}\\ 
e^{im\varphi }&=&\left(-i\text{sgn}\left( \sin \widehat{\varphi }\right) \right) ^{m},
\notag
\end{eqnarray}
and Eq. (\ref{fft12}) takes the form 
\begin{equation}
\!\!\!\!\!\!\!\!\!\!\!\!\!g_{n}^{m}\left( \widehat{\varphi };\frac{\pi }{2}\right)\! =\!\!\sqrt{%
\frac{2n+1}{4\pi }\frac{(n-m)!}{(n+m)!}}\left[i\text{sgn}\left( \sin 
\widehat{\varphi }\right) \right] ^{m}\!\left(\! 1-\frac{i\gamma _{n}^{m}m}{\sin \widehat{\varphi }}\right)\! P_{n}^{m}(\cos \widehat{\varphi }). \label{fft15}
\end{equation}
Note that this expression has a removable singularity at $\sin \widehat{%
\varphi }=0$. Indeed for $m\geqslant 2$ we have $P_{n}^{m}(\cos \widehat{%
\varphi })\sim \sin ^{m}\widehat{\varphi }$, while for $m=1$ we have
\begin{equation}
\left. \frac{P_{n}^{1}(x)}{\sqrt{1-x^{2}}}\right| _{x\rightarrow \pm 1}=-%
\frac{dP_{n}}{dx}\left( \pm 1\right) =-n(n+1)\epsilon _{\pm n},
\label{fft16}
\end{equation}
where symbol $\epsilon _{m}$ is defined by Eq. (\ref{r6}). So, 
\begin{equation}
g_{n}^{m}\left( \pi k;\frac{\pi }{2}\right) =\left\{ 
\begin{array}{c}
\gamma _{n}^{1}\left( -1\right) ^{k+1}\frac{1}{2}\sqrt{\frac{2n+1}{4\pi }%
n(n+1)},\quad m=1, \\ 
0,\quad \quad m\geqslant 2.%
\end{array}
\right.  \label{fft17}
\end{equation}
Hence, the modified algorithm is based on the equation 
\begin{equation}
g_{n}^{m}\left( \widehat{\varphi };\beta \right) =\sum_{m^{\prime
}=-n}^{n}G_{n}^{m^{\prime }m}\left( \beta \right) e^{im^{\prime }\widehat{%
\varphi }},  \label{fft18}
\end{equation}
where function $g_{n}^{m}\left( \widehat{\varphi };\beta \right) $can be
computed for equispaced values of $\widehat{\varphi }$ sampling the full
period. The FFT produces coefficients\ $G_{n}^{m^{\prime }m}\left( \beta
\right) $ from which $H_{n}^{m^{\prime }m}\left( \beta \right) $ can be
found using Eq. (\ref{fft5}).

\section{Numerical experiments}
The algorithms were implemented in Matlab and  tested for $n=0,...,10000$. 

\subsection{Test for recursion stability}

First we conducted numerical tests of the algorithm stability.  Note that the dependence on $\beta $ comes only through the initial
values, which are values of coefficients for layers $m^{\prime }=0$ and $%
m^{\prime }=1$. Hence, if instead of these values we put some arbitrary
function (noise) then we can measure the growth of the magnitude of this noise as the
recursive algorithm is completed. For stable algorithms it is expected that the noise
will not amplify, while amplification of the noise can be measured and some
conclusions about practical value of the algorithm can be made. Formally the
amplitude of the noise can be arbitrary (due to the linearity of
recursions), however, to reduce the influence of roundoff errors we selected
it to be of the order of unity.

\begin{figure}[htb]
\center  \includegraphics[width=0.9\textwidth]{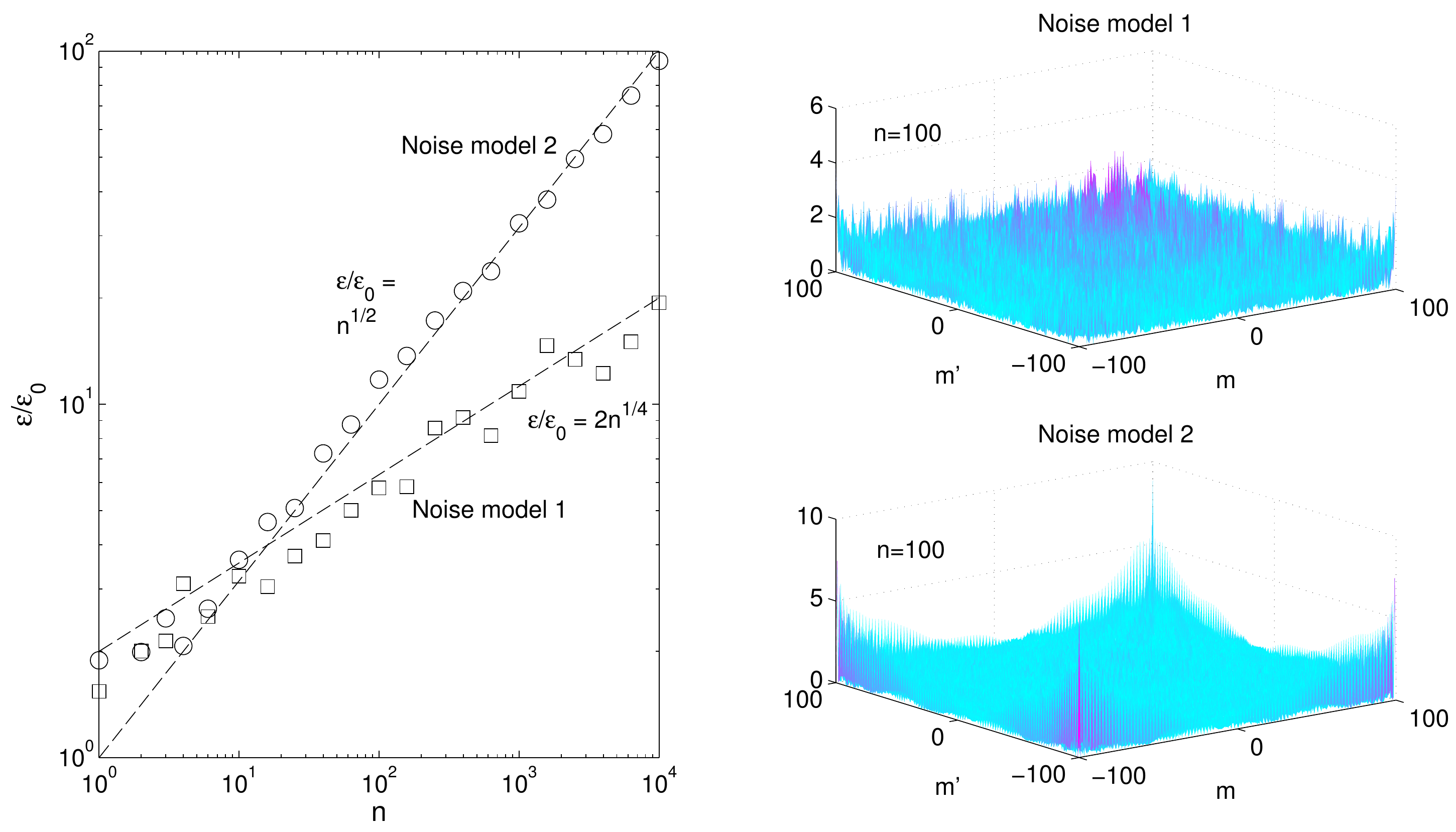}
\caption{The chart on the left shows amplification of the noise in the
proposed recursive algorithms for two noise models. The charts on the right
show distributions of the noise amplitude for some random realization at $%
n=100$.}
\label{Fig6}
\end{figure}

Two models of noise were selected for the test. In the first model
perturbations $\eta _{n}^{m^{\prime },m}$ at the layers $m^{\prime }=0$ and $%
m^{\prime }=1$ were specified as random numbers distributed uniformly
between $-1$ and $1$. In the second model perturbations were selected more coherently. Namely, at $m^{\prime }=0$ the random numbers were non-negative (distributed between $0$ and $1$). At $m^{\prime }=1$ such a random
distribution was pointwise multiplied by factor $\left( -1\right) ^{m}$. The reason for this factor is that effectively this brings some symmetry for resulting distributions of $\eta _{n}^{m^{\prime },m}$ for $m^{\prime }>0$
and $m^{\prime }<0$. Figure \ref{Fig6} shows that in the first noise model the overall error (in the $L^{\infty }$ norm) grows as $\sim n^{1/4}$, while
for the second noise model the numerical data at large enough $n$ are well approximated by $\epsilon =\epsilon _{0}n^{1/2}$. Note that the data points shown on this figure were obtained by taking the maximum of 10 random realizations per each data point. On the right hand side of Fig. \ref%
{Fig6} are shown error distributions for some random realization and some $n$ ($n=100$, the qualitative picture does not depend on $n$). It is seen that for
the first noise model the magnitude of $\eta _{n}^{m^{\prime },m}$ is
distributed approximately evenly (with higher values in the central region and diagonals $m^{\prime }=\pm m)$. For the second noise model the distribution is substantially different. The highest values are observed in the boundary regions $m^{\prime }\approx \pm n$ and $m\approx \pm n$ with
the highest amplitudes near the corners of the computational square in the $\left( m,m^{\prime }\right) $ space. 

The behavior observed in the second noise model can be anticipated, as the
scheme is formally unstable, the absolute values of coefficients 
$a$ and $b$ (see Eq. (\ref{Neu1.1})) grow near the boundaries of the computational domain and a fast change of these coefficients near the boundaries requires some other technique for investigation of instabilities than the method used. Smaller errors and their distribution observed in the first noise model are more puzzling, and  we can
speculate about some cancellation effects for random quantities with zero mean
appearing near the boundaries, and to the variability of coefficients $%
a,b,$ and $c$ in Eq. (\ref{Neu1.1}), so that the stability analysis is only approximate. What is important that in all our tests with different
distributions of initial values of $\eta _{n}^{m^{\prime },m}$ we never
observed exponential growth. The maximum growth rate behaved at large $n$ as $%
n^{\alpha }$, $\alpha \approx 1/2$. Hence, for $n\sim 10^{4}$ one can expect
the errors in the domain of two orders of magnitude larger than the errors
in the initial conditions, which makes the algorithm practical. Indeed, in double precision, which provides errors $\sim 10^{-15}$
in the initial values of the recursions, then for $n=10^{4}$ one can expect
errors $\sim 10^{-13}$, which is acceptable for many practical problems. Of
course, if desired the level of the error can be reduced, if needed, using e.g. quadruple precision, etc.

\subsection{Error and performance tests}

The next error tests were performed for actual computations of $%
H_{n}^{m^{\prime },m}$. For small enough $n$ ($n\sim 10$) one can use an
exact expression (\ref{r5}) as an alternative method to figure out the
errors of the present algorithm. Such tests were performed and absolute
errors of the order of $10^{-15}$, which are consistent with double
precision roundoff errors were observed. The problem with sum (\ref{r5}) is
that at large $n$ it requires computation of factorials of large numbers,
which creates numerical difficulties. While computation of factorials and
their summation when the terms have the same sign is not so difficult (e.g.
using controlled accuracy asymptotic expansion), the problem appears in the
sums with large positive and negative terms. In this case to avoid the loss
of information special techniques of working with large integers (say, with
thousand digits) should be employed. This goes beyond the present study, and
we used different methods for validation than comparing with these values.

Another way is to compute $H_{n}^{m^{\prime },m}\left( \beta \right) $ is based on the flip decomposition, i.e. to use one of equations (\ref{af5})-(\ref{af7}). In this decomposition all coefficients are ``good'' in terms that complex exponents or cosines can be computed accurately. These formulae
require the flip rotation coefficients, $H_{n}^{m^{\prime }m}\left( \pi /2\right) $, which should be computed and stored
to get $H_{n}^{m^{\prime },m}\left( \beta \right) $ for arbitrary $\beta $.
Note that these relations hold also for $\beta =\pi /2$,
which provides a self-consistency test for $H_{n}^{m^{\prime }m}\left( \pi
/2\right) $. Despite the summations in (\ref{af5})-(\ref{af7}) requires $O(n^{3})$
operations per subspace $n$ and are much slower than the algorithms proposed
above, we compared the results obtained using the recursive algorithm for
consistency with Eqs. (\ref{af5}) and found a good agreement (up to the numerical errors reported below) for $n$ up to 5000 and different $\beta $
including $\beta =\pi /2$. 

One more test was used to validate computations,
which involve both recursions and relation (\ref{af7}). This algorithm with
complexity $O\left( n^{2}\right) $ per subspace was proposed and tested in \cite{Gumerov12:JCC}. There coefficients $%
H_{n}^{m^{\prime }m}\left( \pi /2\right) $ were found using the present
recursion scheme, which then were used only to compute the diagonal
coefficients $H_{n}^{mm}\left( \beta \right) $ and $H_{n}^{m,m+1}\left(
\beta \right) $ at arbitrary $\beta $. The recursion then was applied to
obtain all other coefficients, but with propagation from the diagonal
values, not from $H_{n}^{0,m}\left( \beta \right) $ and $H_{n}^{1,m}\left(
\beta \right) $. Motivation for this was heuristic, based on the
observation that $O\left( 1\right) $ magnitudes are achieved on the diagonals and then they can decay exponentially (see Fig. \ref{Fig2}), so it
is expected that the errors will also decay. The tests for $n$ up to 5000
showed that such a complication of the algorithm is not necessary, and both
the present and the cited algorithms provide approximately the same errors.

In the present study as we have two alternative ways to compute\ $%
H_{n}^{m^{\prime }m}\left( \beta \right) $ using the recursive algorithm and the
FFT-based algorithm, we can use self-consistency and cross validation tests
to estimate the errors of both methods.

Self consistency tests can be based on relation (\ref{b1}). In this case
computed $H_{n}^{m^{\prime }m}\left( \beta \right) $ were used to estimate
the following error 
\begin{equation}
\epsilon _{n}^{(0)}\left( \beta \right) =\max_{m,m^{\prime }}\left|
\sum_{\nu =-n}^{n}H_{n}^{m^{\prime }\nu }\left( \beta \right) H_{n}^{\nu
m}\left( \beta \right) -\delta _{mm^{\prime }}\right| ,\quad n=0,1,...
\label{ca5}
\end{equation}
We found also the maximum of $\epsilon _{n}^{(0)}\left( \beta \right) $ over
five values of $\beta =0,\pi /4,\pi /2,3\pi /4,\pi $ for the recursive
algorithm and for two versions of the FFT-based algorithms. For the basic
FFT-based algorithm we used $\widehat{\theta }=\pi /2-\xi /n$, where $\xi $
was some random number between 0 and 1. For the modified algorithm, which
has some arbitrary coefficient $\gamma _{n}^{m}$ we used $\gamma
_{n}^{m}=1/n $, which, as we found provides smaller errors than $\gamma
_{n}^{m}=1$ or $\gamma _{n}^{m}=1/n^{2}$ and consistent with the
consideration of magnitude of the odd and even normalization coefficients
(see Eq. (\ref{fft7})). The results of this test are presented in Fig. \ref%
{Fig7}. It is seen that while at small $n$ the error is of the order of the
double precision roundoff error at larger $n$ it grow as some power of $n$.
The error growth rate at large $n$ for the recursive algorithm is smaller
and approximately $\epsilon _{n}^{(0)}\sim n^{1/2}$, which is in a good
agreement with the error growth in the noise model \#2 discussed above. For
the FFT-based algorithms the error grows approximately as $\epsilon
_{n}^{(0)}\sim n^{3/2}$, so it can be orders of magnitude larger than in the
recursive algorithm, while still acceptable for some practical purposes up
to $n\sim 10^{3}$. Such grows can be related to summation of coefficients of
different magnitude in the FFT, which results in the loss of information.
Comparison of the basic and the modified FFT-based algorithms show that the
errors are approximately the same for the both versions, while the error in
the basic algorithm can behave more irregularly than that in the modified
algorithm. This can be related to the fact that used values of $\widehat{%
\theta }$ at some $n$ were close to zeros of the associated Legendre
functions, and if this algorithm should be selected for some practical use
then more regular way for selection of $\widehat{\theta }$ should be worked
out. 
\begin{figure}[htb]
\center  \includegraphics[width=0.9\textwidth]{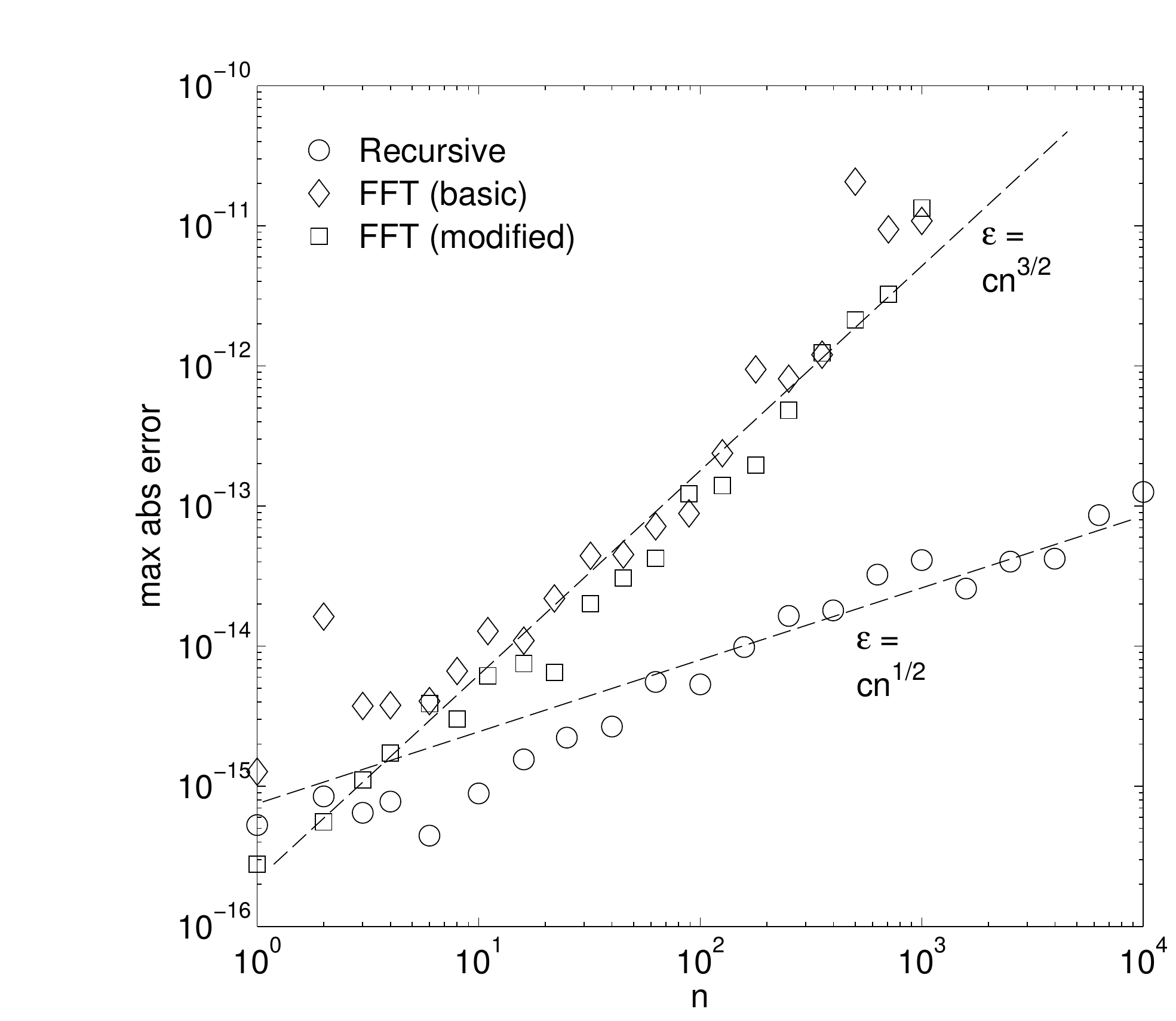}
\caption{Self-consistency error test of the recursive and FFT-based
algorithms validating that the symmetric matrix of rotation coefficients is
unitary.}
\label{Fig7}
\end{figure}

\begin{figure}[htb]
\center  \includegraphics[width=0.9\textwidth, trim =2in 0  0.5in 0in]{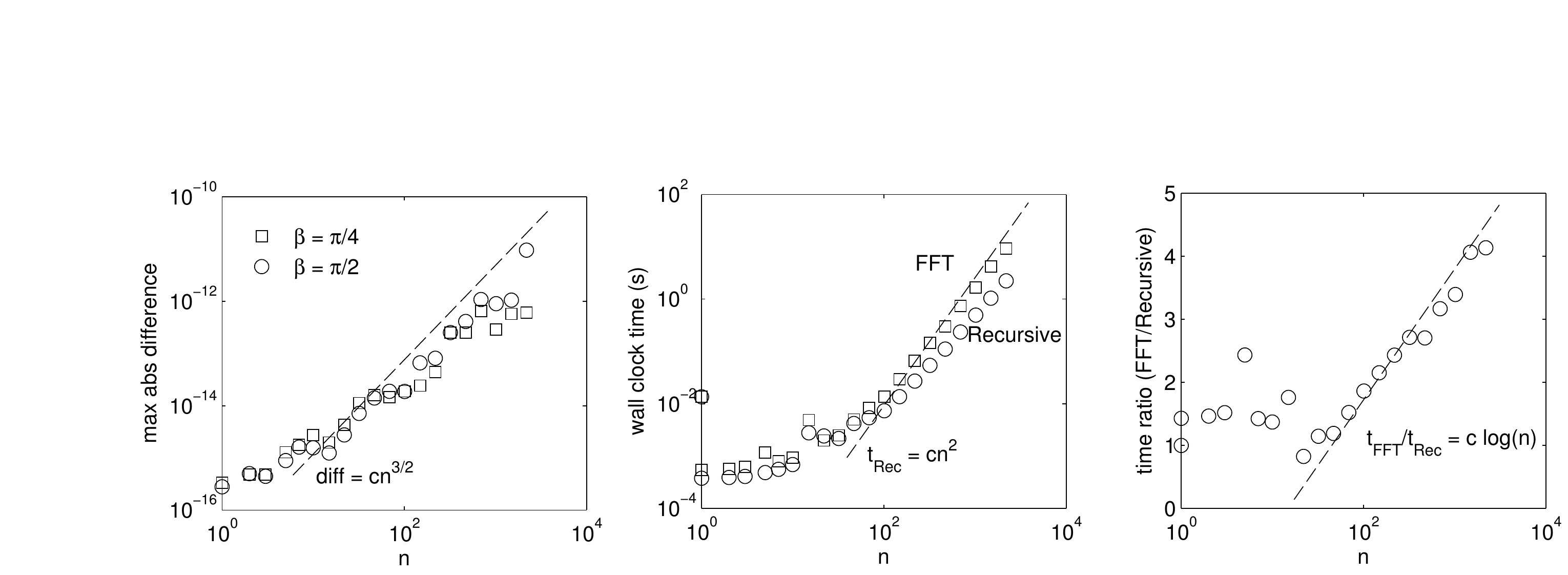}
\caption{The maximum of absolute difference in the rotation coefficients $%
H_{n}^{m^{\prime }m}$ computed using the recursive and the FFT-based
(modified) algorithms as a function of $n$ at two values of $\protect\beta $
(left). The center and the right plots show 
wall-clock times for these algorithms and the ratio of these times, respectively (standard PC, Matlab).}
\label{Fig8}
\end{figure}

The second test we performed is a cross-validation test. In this case we
computed the difference 
\begin{equation}
\epsilon _{n}^{(1)}\left( \beta \right) =\max_{m,m^{\prime }}\left|
H_{n}^{m^{\prime }m(FFT)}\left( \beta \right) -H_{n}^{m^{\prime
}m(rec)}\left( \beta \right) \right| ,\quad n=0,1,...  \label{ca6}
\end{equation}
We also measured and compared the wall clock times for execution of the
algorithms (standard Matlab and its FFT library on a standard personal
computer). Results of these tests are presented in Fig. \ref{Fig8}. First we
note that both FFT-based algorithms showed large errors for $n$ $\symbol{126}
$ $2300$ and failed to produce results for larger $n$. This can be related
to the loss of information in summation of terms of different magnitude, as
it was mentioned above. On the other hand the recursive algorithm was
producing reasonable results up to $n=10000$ and there were no indications
that it may not run for larger $n$ (our constraint was the memory available
on the PC used for the tests). So the tests presented in the figure were
performed for $n\leqslant 2200.$ It is seen that the difference is small
(which cross-validates the results in this range), while $\epsilon
_{n}^{(1)}\left( \beta \right) $ grows approximately at the same rate as the
error $\epsilon _{n}^{(0)}$ for the FFT-based algorithms shown in Fig. \ref%
{Fig7}. Taking into account this fact and numerical instability of the
FFT-based algorithms we relate it rather to the errors in that algorithms,
not in the recursive algorithm. We also noticed that the FFT-based algorithm
at large $n$ produces somehow larger errors for $\beta =\pi /2$ than for
other values tried in the tests. In terms of performance, it is seen that
both, the recursive and the FFT-based algorithm are well-scaled and $O\left(
n^{2}\right) $ scaling for large enough $n$ is acieved by the recursive
algorithm nicely. The FFT-based algorithm in the range of tested $n$ can be
several times slower. The ratio of the wall-clock times at $n>100$ is well
approximated by a straight line in semi-logarithmic plots, which indicates $%
\log n$ behavior of this quantity, as expected. 

\section{Conclusion}
This paper first presented a study of the asymptotic behavior of the rotation coefficients $%
H_{n}^{m^{\prime }m}\left( \beta \right) $ for large degrees $n$. Based on
this study, we proposed a recursive algorithm for computation of these
coefficients, which can be applied independently for each subspace $n$
(with cost $O\left( n^{2}\right) $) or to $p$ subspaces ($n=0,...,p-1$)
(cost $O\left( p^{3}\right) $). A theoretical and numerical analysis
of the stability of the algorithm shows that while the algorithm is weakly unstable, the growth rate of perturbations is small enough, which makes it
practical for computations for relatively large $n$ (the tests were performed up to $n=10^{4}$, but the scaling of the error indicates that even larger $n$ can be computed). Alternative FFT-based algorithms of complexity $O\left(
n^{2}\log n\right) $ per subspace $n$ were also developed and studied. Both the FFT-based and recursive computations produce consistent results for $n\lesssim 10^{3}$. In this range the recursive algorithm is faster and produces smaller
errors than the FFT-based algorithm.

\end{document}